\documentclass{article}

     \usepackage[preprint,nonatbib]{neurips_2020}

\usepackage[utf8]{inputenc} 
\usepackage[T1]{fontenc}    
\usepackage{hyperref}       
\usepackage{url}            
\usepackage{booktabs}       
\usepackage{amsfonts}       
\usepackage{nicefrac}       
\usepackage{microtype}      

\usepackage{graphicx}
\usepackage{subfigure}
\usepackage{booktabs} 
\usepackage{amssymb}
\usepackage{amsmath}
\usepackage[ruled, vlined, linesnumbered]{algorithm2e}
\usepackage{algorithmic}
\usepackage{multirow}
\newtheorem{theorem}{Theorem}[section]

\newcommand*{\QEDB}{\hfill\ensuremath{\square}}%
\newcommand\vect[1]{{#1}}
\newcommand{\dimx}{\mathbb{R}^{M}}
\newcommand{\dimy}{\mathbb{R}^{N}}
\newcommand{\dimalpha}{\mathbb{R}^{m}}
\newcommand{\dimbeta}{\mathbb{R}^{n}}
\newcommand\Item[1][]{%
  \ifx\relax#1\relax  \item \else \item[#1] \fi
  \abovedisplayskip=0pt\abovedisplayshortskip=0pt~\vspace*{-\baselineskip}}
  
\DeclareMathOperator*{\argmin}{arg\,min}
\title{Primal-Dual Sequential Subspace Optimization for Saddle-point Problems}

\author{%
    Yoni Choukroun \\
    Huawei\\
    \texttt{choukroun.yoni@gmail.com}
    \And
    Michael Zibulevsky\\
    Technion, IIT\\
    \texttt{mzibul@gmail.com}
    \And
    Pavel Kisilev\\
    Huawei\\
    \texttt{pavel.kisilev@huawei.com}
}

\begin{document}

\maketitle

\begin{abstract}
We introduce a new sequential subspace optimization method for large-scale saddle-point problems.
It solves iteratively a sequence of auxiliary saddle-point problems in low-dimensional subspaces, spanned by directions derived from first-order information over the primal \emph{and} dual variables.
Proximal regularization is further deployed to stabilize the optimization process.
Experimental results demonstrate significantly better convergence  relative to popular first-order methods. 
We analyze the influence of the subspace on the convergence of the algorithm, and assess its performance in various deterministic optimization scenarios, such as bi-linear games, ADMM-based constrained optimization and generative adversarial networks.
\end{abstract}

\section{Introduction}
\label{Introduction}
Saddle-point problems arise in many applications, such as game theory \cite{leyton2008essentials}, constrained and robust optimization \cite{arrow1958studies,ben2009robust} and generative adversarial networks (GANs) \cite{goodfellow2014generative}.
Important variational problems such as $\ell_{\infty}$ minimization, convex segmentation or compressed sensing \cite{chan2006algorithms,candes2006stable} have saddle-point formulations that are efficiently handled using primal-dual solvers.

Convex-concave saddle point problems have been widely investigated in the context of constrained optimization problems that are reduced to finding the saddle-point of a corresponding Lagrangian.
Gradient methods are computationally efficient and proven to be convergent under strict convexity or concavity of the primal-dual variables \cite{arrow1958studies,popov1980modification,chambolle2011first,cherukuri2017saddle}.
Nevertheless, there is a broad spectrum of problems wherein a function considered is not strictly convex-concave. It has been observed that when using the naive gradient method, convergence is not guaranteed, leading instead to diverging oscillatory solutions \cite{arrow1958studies,nedic2009subgradient,feijer2010stability,ImproveTrainGans}.
Following the solution of \cite{nesterov2005smooth} on specific convex-concave setting, \cite{nemirovski2004prox} proposed a proximal method \cite{rockafellar1976monotone} for saddle-point problem with sub-linear convergence rate. 
Later, \cite{nesterov2009primaldual} provided a worst-case optimal dual averaging scheme with proven sub-linear convergence for general smooth convex-concave problems.
In the emerging topic of generative adversarial networks, many new approaches have been proposed or re-investigated \cite{korpelevich1976extragradient,kulis2010implicit,rakhlin2013optimization,daskalakis2017training,mescheder2017numerics}.
However, gradient based methods suffer from the coupling term usually present in min-max games \cite{liang2018interaction}, and remain inherently slow especially in ill-conditioned problems.

In case of large scale optimization problems, there is a need for optimization algorithms whose storage requirement and computational cost per iteration grow at most linearly with the problem dimensions. 
In the context of minimization, this constraint has led to the development of a broad family of methods such as variable metric methods, and \emph{subspace optimization}. 
The most simple (and ubiquitous) subspace minimization algorithm is steepest descent coupled with line-search procedure. Early methods proposed to extend the minimization to a $k$ dimensional subspace spanned by various directions, such as gradients, conjugate directions, previous outer iterations or Newton directions, see \cite{cragg1969study,miele1969study,dennis1987generalized,conn1996iterated}.
One popular algorithm of this type is the Conjugate Gradient (CG) method \cite{hestenes1952methods} that possesses a remarkable linear convergence rate induced by the well-known expanding manifold property of \emph{quadratic} CG.
In the case of general smooth convex functions, preserving the expanding manifold property enables linear convergence rate of CG, but the cost of each iteration increases with the iteration count, and thus makes it computationally prohibitive.
To alleviate this problem, \cite{nemirovski1982orth} suggested to restrict the optimization subspace to the span of the current gradient, the average of previous steps and a weighted sum of previous gradients.
The resulting ORTH-method inherits the optimal worst-case convergence, but fails to preserve CG convergence rate for quadratic objective, leading to slower convergence in the neighborhood of a solution.
To address this problem the SEquential Subspace OPtimization (SESOP) method \cite{Narkiss-2005} extends the ORTH subspaces with the directions of the last propagation steps. 


In this work, motivated by the inherent slowness of gradient based methods and the power of subspace optimization, we extend the idea of subspace optimization, until now limited to minimization, to saddle-point problems. 
Specifically, we solve sequentially low dimensional saddle-point problems in subspaces defined by first-order information. 
We propose to perform the subspace optimization over the primal \emph{and} dual variables, allowing to search for a saddle-point in a richer subspace, wherein the function can increase and/or decrease in primal and dual variables respectively. 
Further, we propose to couple the saddle-point objective with proximal operators in order to ensure the existence of a stationary point in the subspace. 
We solve the subspace optimization via adapted second order optimization that can be implemented efficiently in the given low dimensional subspace.
Finally, we perform backtracking line search over the gradient norm. This ensures faster convergence, and most importantly, prevents divergence in degenerative cases.
Experimental results assess the power and usefulness of the proposed method.
\section{Background}
\label{Background}
As our approach adapts subspace optimization to saddle point problems, we start with notations related to the saddle-point setting, and provide brief review of subspace methods for minimization.
\subsection{Notation and Definitions}
We consider the unconstrained saddle-point problem
\begin{equation} \label{eq:saddle_point_formulation}
\min_{\vect{x} \in \dimx}\max_{\vect{y} \in \dimy} f(\vect{x,y}),
\end{equation}
where $f$ is twice continuously differentiable  and the first derivative is $L$ Lipschitz continuous on $\dimx \times \dimy$.
The point $(\vect{x^*,y^*})$ is a global saddle-point iff
\begin{equation} \label{eq:saddle_point_optimality}
f(\vect{x^*,y})\leq f(\vect{x^*,y^*})\leq f(\vect{x,y^*})\ \ \ \forall \vect{x}\in \dimx, \vect{y}\in \dimy.
\end{equation}
Finding a global saddle-point is computationally intractable. We therefore assume local convexity-concavity of the objective in which case eq. (\ref{eq:saddle_point_optimality}) holds in a \emph{local} $r$ neighborhood ball of $\vect{(\vect{x^*,y^*})}$, which we denote $B_{2}((\vect{x^*,y^*}),r)$.
We define $(x^* ,y^*)$ to be a local \emph{stable} saddle-point when the gradient is zero and the second order derivatives are positive definite in $x$ and negative definite in $y$.
This notion of stability is stronger than in regular saddle-point (i.e. second derivatives are positive and negative \emph{semi-} definite matrices respectively), since the maximum and minimum eigenvalue of $\nabla_{xx} f(x^{*},y^{*})$ and $\nabla_{yy} f(x^{*},y^{*})$ respectively are bounded away from 0.
We further define $(x^* ,y^*)$ to be a local \emph{unstable} saddle-point when the gradient is zero and the second order derivatives in $x$ and $y$ are vanishing.
\subsection{Subspace Optimization for Minimization}
The core idea of subspace optimization is to perform the optimization of a function in a small subspace spanned by a set of directions obtained from an available oracle.
Denoting a function $g:\mathbb{R}^n\rightarrow\mathbb{R}$ to be minimized and $P_k\in \mathbb{R}^n\times\mathbb{R}^d, d<< n$ as the set of $d$ directions at iteration $k$, an iterated subspace optimization method aims at solving the following minimization problem
\begin{equation}\label{eq:sesop_minimization}
\begin{aligned}
	\alpha_k = \argmin_{\alpha \in \mathbb{R}^d} g(x_k + P_k \alpha),
\end{aligned}
\end{equation}
followed by the update rule $x_{k+1} = x_k +P_k \alpha_k$.
The dimensions of the problem are then reduced from the optimization space  $\mathbb{R}^n$, to the controlled subspace in $\mathbb{R}^d$ spanned by the columns of $P_k$. 

The subspace structure may vary depending on the chosen optimization technique.
Krylov descent method defines the subspace as $\text{span}\{H^{0}\nabla g,...,H^{d-1}\nabla g\}$ for some preconditioning matrix $H$, e.g. $H=\nabla^2g$ in \cite{vinyals2012krylov}.
Related to Krylov subspaces, the Conjugate Gradient method \cite{hestenes1952methods} reduces the search space to current gradient and previous step, such that $\text{span}\{p_{k},\nabla g\}$, where $p_{k}=x_{k}-x_{k-1}$.
Nemirovski \cite{nemirovski1982orth} provided optimal worst case complexity of the method in convex setting by defining the subspace as $\text{span}\{\nabla g(x),x_{k}-x_{0},\sum_{j=0}^{k}w_{j}\nabla g(x_{j})\}$,
with appropriate weights $\{w_{j}\}_{j=0}^{k}$.
The SESOP algorithm \cite{Narkiss-2005} extends this method by adding the previous search directions $\{p_{k-i}\}_{i=0}^{d-3}$. Thus, the method provides optimal worst case complexity, generalizes the quadratic CG method, and allows truncated approximation of the expanding manifold property on non-linear objectives.
Also, more directions can be added to enrich the subspace in order to improve the convergence \cite{conn1996iterated,ZibEladSPM}. 

\section{Sequential Subspace Optimization for saddle-point Problems}
We extend here the iterative formulation of the subspace optimization to saddle-point problems.
We first formulate then motivate the use of primal-dual subspace optimization, and then provide adapted generic optimization framework.
\subsection{Formulation}
Let us define the subspace saddle-point problem as 
\begin{equation}\label{eq:sesop_minimax}
\begin{aligned}
\min_{\vect{\alpha} \in \dimalpha}\max_{\vect{\beta} \in \dimbeta}f(x_{k}+P_k\alpha,y_{k}+Q_k\beta),
\end{aligned}
\end{equation}
where we assume $m \ll M$ and $n \ll N$. Matrices $P_k$ and $Q_k$ define the subspace structure at iteration $k$.
The subspace optimization can be solved exactly \emph{or} approximately.
The new iterate is of the form
\begin{equation}\label{eq:sesop_update}
\begin{aligned}
                \begin{pmatrix}
                  x_{k+1}  \\
                    y_{k+1} 
                \end{pmatrix}
                =
                \begin{pmatrix}
                  x_{k} \\
                  y_{k}
                \end{pmatrix}
                +\eta_k
                \begin{pmatrix}
                  P_k\alpha\\
                    Q_k\beta
                \end{pmatrix}
                ,
\end{aligned}
\end{equation}
where $\eta_{k}$ is the step size obtained via outer optimization (i.e. original problem), and the procedure stops if convergence tests are satisfied.
This formulation allows flexibility in definition of the search space. The following simple but challenging example illustrates this property. The  bi-linear game $f(x,y)=x^{T}y$ \cite{ImproveTrainGans} diverges when search is performed over one dimensional anti-gradient/gradient direction. However, convergence can be reached if the optimization is performed separately over the primal and dual variables, as shown in the following theorem.
Proofs are provided in the Appendix.
\begin{theorem}\label{theorem:xy}
Consider the saddle-point problem $f(x,y)=x^{T}y$ and the update from eq. (\ref{eq:sesop_update}),
where $\alpha$ and $\beta$ are obtained by solving eq.(\ref{eq:sesop_minimax}) with $P_k=\nabla_x f(x_k,y_k)$ and  $Q_k=\nabla_y f(x_k,y_k)$.
Then, $\forall \eta_k \in (0,{2f(x_k,y_k)^{2}}/{\|P_k\|^{2}\|Q_k\|^{2}})$ the procedure converges to optimum.\\
Also, the gradient method (i.e. $-\alpha=\beta>0$), diverges $\forall \eta_k>0$. 
\end{theorem}
We emphasize the fact that \emph{joint} subspace optimization is convergent while independent (alternating) subspace optimization is divergent in this unstable case.
Following the subspace minimization strategy to use more than current gradient, we seek a saddle-point in the subspace spanned by first order information. Namely, we use the \emph{mandatory} \cite{nemirovski1982orth} current gradient, previous gradients, and the previous search steps in $\vect{x}$ and $\vect{y}$, such that $\text{span}\{P_{k}\}= \text{span}\{S^{x}_{k},G^{x}_{k}\}$ and $\text{span}\{Q_{k}\}= \text{span}\{S^{y}_{k},G^{y}_{k}\}$,
where $S^{u}_{k}=\{p^{u}_{k-l-1},\dots,p^{u}_{k-1}\}$ and $G^{u}_{k}=\{\nabla_u f(x_{k-l},y_{k-l}),\dots,\nabla_u f(x,y)\}$, with $p^{u}_{k}=u_{k}-u_{k-1}$.
Other directions can be used or added to improve the convergence as well.
Expanding the subspace with more directions can enrich the subspace but enables a subjective trade-off between computational cost and speed of convergence.
Such subspace formulation generalizes popular methods, e.g. the gradient method \cite{rockafellar1976monotone} or Optimistic Mirror Descent \cite{rakhlin2013optimization,daskalakis2018limit}. In order to improve the convergence, the proposed framework can be combined with other methods, such as weighted averaging of iterates as final solution \cite{bruck1977weak}, or consensus optimization \cite{mescheder2017numerics} as modification of the objective.

\subsection{General Subspace Convergence Analysis}\label{subsec:conditions}
In this section we analyse the convergence conditions of the subspace optimization method, in terms of the norm of gradient, i.e. convergence to stationary point.
Consider $z_{k}=[x_{k},y_{k}]$
where $[\cdot,\cdot]$ denotes vectors concatenation.
We denote the direction $d_{k}=[P_{k}\vect{\alpha,Q_{k}\beta}]=R_{k}[\alpha,\beta]:=R_{k}\gamma$, where $R_{k}$ is the block matrix populated with $P_{k}$ and $Q_{k}$ in the block diagonal and zero elsewhere.
Here we assume $R_k$ has linearly independent columns. 
From the first order expansion, there exists sufficiently small $\zeta>0$ such that
\begin{equation}\label{eq:taylor}
\begin{aligned}
f(z_{k+1}) &= f(z_{k}+\zeta d_{k})= f(z_k) + \zeta \langle\nabla f(z_k),d_k\rangle + o(\zeta^{2}\|d_{k}\|).
\end{aligned}
\end{equation}
Thus, by taking derivative of eq.(\ref{eq:taylor}) we have $\nabla{f}(z_{k+1}) \approx \nabla f(z_k) + \zeta \nabla^{2}f(z_k)d_k$.
Since we are interested in decreasing the gradient norm to reach convergence, we have
\begin{equation}\label{eq:conv}
\begin{aligned}
&\|\nabla f(z_{k+1})\|^2 =  \|\nabla f(z_k)\|^2 + 2\zeta \nabla f(z_k)^{T}\nabla^{2} f(z_k)d_k +o(\zeta^{2}\|\nabla^{2} f(z_{k})d_k\|).
\end{aligned}
\end{equation}
Thereafter, the sufficient condition for convergence to local stationary point is $\nabla f(z_k)^{T}\nabla^{2} f(z_k)d_k<0$.
For example, the steepest descent/ascent direction is convergent in the case of strongly convex-concave problem since in that case $\nabla f(z_k)^{T}\nabla^{2} f(z_k)\nabla f(z_k)<0$.
We can reformulate the previous equation in terms of the subspace parameters such that, 
since $\nabla_{\gamma}f(z+R\gamma)=R^{T}\nabla_{z}f(z+R\gamma)$, we have
\begin{equation}\label{eq:conv1}
\begin{aligned}
&\|\nabla f(z_{k+1})\|^2 \approx \|\nabla f(z_{k})\|^2 +  2\zeta \nabla f(z_{k})^{T}\nabla^{2} f(z_{k})R_{k}\gamma\\
=&\|\nabla f(z_{k})\|^2 + 2\zeta \nabla_{\gamma} f(z_{k}+R_{k}\gamma)^{T}\big\rvert_{\gamma=0}R^{+T}_{k}R^{+}_{k}\nabla^{2}_{\gamma} f(z_{k}+R_{k}\gamma)\big\rvert_{\gamma=0}\gamma,
\end{aligned}
\end{equation}
where, $R_{k}^{+}$ denote the Moore–Penrose pseudoinverse of matrix $R_{k}^{T}$.
Thus, assuming single Newton step in subspace $\gamma=-\nu\nabla^{2}_{\gamma} f(z_{k}+R_{k}\gamma)^{-1}\big\rvert_{\gamma=0}\nabla_{\gamma} f(z_{k}+R_{k}\gamma)\big\rvert_{\gamma=0}$, $\exists \nu>0$ such that 
\begin{equation}\label{eq:conv2}
\begin{aligned}
\|\nabla f(z_{k+1})\|^2&= \|\nabla f(z_k)\|^2- 2\nu\zeta \|R^{+}_{k}\nabla_{\gamma} f(z_{k}+R_{k}\gamma)\big\rvert_{\gamma=0}\|^2.
\end{aligned}
\end{equation}
According to eq (\ref{eq:conv2}), in the neighborhood of the current point $z_{k}$, Newton step in the subspace domain decreases gradient norm of the original problem, and thus induces global convergence to stationary point. 
However, \emph{contrary to the minimization setup} \cite{conn1996iterated} where subspace optimization are descent methods,
exact convergence to stationary point in subspaces (i.e. $\nabla_{\gamma}f(z_{k}+R_{k}\gamma)=0$) does not necessarily enforce convergence in the original problem space, as shown in Theorem \ref{theorem:xy}.
This is mainly due to the interaction term $\nabla_{xy}f(x,y)$ present in eq. (\ref{eq:conv}) via the Hessian matrix.
Thus, the extension of subspace optimization to saddle-point problems is then not straightforward.
\subsection{Convergence Improvement Strategies}
To ensure convergence of inner and outer optimization, we propose to solve the subspace optimization in a constrained local region. Also, we propose to correct the direction obtained from the subspace optimization by controlling the outer step size $\eta_{k}$ as in eq. (\ref{eq:sesop_update}) via an adapted line-search procedure.
\subsubsection{Proximal Regularization}
In the following, we extend the auxiliary subspace saddle problem (\ref{eq:sesop_minimax}) by adding proximal point regularization \cite{martinet1970breve}.
At each iteration we solve the subspace proximal problem over $\tilde{f}(x,y)$, namely
\begin{equation}\label{eq:sesop_minimax_prox}
\begin{aligned}
\min_{\vect{\alpha} \in \dimalpha}\max_{\vect{\beta} \in \dimbeta}
\tilde{f}(x_{k}+P_k\alpha,y_{k}+Q_k\beta):=&f(x_{k}+P_k\alpha,y_{k}+Q_k\beta) \\+& \frac{\tau_{k}}{2} \|x_{k}+P_k\alpha-\bar{x}_k\|^2- \frac{\tau_{k}}{2} \|y_{k}+Q_k\beta-\bar{y}_k\|^2,\end{aligned}
\end{equation}
where $\bar{x}_k$ and $\bar{y}_k$ denote the primal and dual prox-centers respectively (e.g. moving average or previous point).
The proximal approach motivation is two-fold. First, it allows averaging over iterations, reducing the oscillation behavior typical to min-max games \cite{ImproveTrainGans}, and improves stability of the optimization procedure.
Foremost, it ensures the existence of a saddle-point in potentially degenerate subspaces, avoiding divergence, in a trust-region fashion \cite{sorensen1982newton}. 
\subsubsection{Saddle-point Backtracking Line Search}
Common line-search backtracking methods \cite{Nocedal2006} cannot be applied straightforwardly to the saddle-point problems,
since implementing search over the function primal and dual values can diverge (e.g. bi-linear game).
To tackle this problem, we perform backtracking line search over the gradient norm to both ensure faster convergence of the method, and, most importantly prevent potential divergence after the inner subspace optimization.
The proposed procedure is described in Algorithm \ref{algo:line_search}.
\begin{algorithm}
\caption{Saddle Backtracking Line Search}\label{algo:line_search}
\begin{algorithmic}
\REQUIRE{$c\in [0,1),\nu\in (0,1), \eta\leq 1$}
\STATE {\bfseries Input:} $f:\mathbb{R}^{M+N}\rightarrow \mathbb{R}$, current point $z_{k}$, direction $d_k$
\STATE {\bfseries Output:} Step size $\eta$
\WHILE{$\|\nabla {f}(z_{k}+\eta d_k)\|^2\geq\|\nabla f(z_k)\|^{2}+\eta c \nabla f(z_k)^{T}\nabla^{2} f(z_k)d_k$}
\STATE $\eta = \eta*\nu$
\ENDWHILE
\STATE \textbf{return} $\eta$
\end{algorithmic}
\end{algorithm}
This step size search procedure is used in \emph{both} inner (fast convergence in subspace) and outer (step correction) optimization.
In our experiments we chose $c=0$ for less computational overhead and set $\nu=0.5$. We further limit the number of line-search iterations to 30.
The following theorem, based on the analysis of Section \ref{subsec:conditions}, states the convergence of the proposed algorithm for the standard gradient method (memoryless subspaces).
\begin{theorem}\label{theorem:linesearch}
Consider function $f(x,y)$ with \emph{stable} saddle-point $(x^*,y^*)$.
Assume the subspace is spanned by the anti-gradient and gradient directions for the primal and dual variables, respectively.
Then, the procedure of Algorithm \ref{algo:line_search} converges to the optimum for every $(x_0,y_0)\in B_2((x^*,y^*),r)$.
\end{theorem}
\subsection{Second-Order Saddle-point Optimization in Subspace}
Second order methods aim at finding roots of the gradient via solution of the second order expansion.
Therefore, they can converge extremely fast to saddle-points, especially in the proximity of the solution where the problem has a good quadratic approximation.
The major drawback is the prohibitive computational cost for both the computation and inversion of the Hessian.
However, in our small subspace setting, second order methods, s.a. (Quasi-)Newton, can be handled efficiently. 
In particular, the computation of the Hessian in the subspace is performed via Hessian product with the direction vectors \cite{Pearlmutter94fastexact}. 
Nowadays, it can be handled efficiently via automatic differentiation tool, since $(\partial^{2}f)\cdot v=\partial(\partial f \cdot v)$. 
Hessian inversion is computationally negligible in \emph{low dimensional} subspace (generally up to ten dimensions).
The method can be further accelerated using frozen or truncated Hessian strategies, especially when the Hessian remains almost unchanged in the vicinity of the solution.
The second order proximal subspace optimization is performed iteratively until the convergence (or maximum number of iterations) is reached, as follows:
\begin{equation}\label{eq:sesop_minimax_tau}
\begin{aligned}
&\gamma_{k+1} = \gamma_k - \eta_k \tilde{H}_{\gamma}^{-1}(z_{k}+R_{k}\gamma_k) \nabla_{\gamma} \tilde{f}(z_{k}+R_{k}\gamma_k)\\
			 &= \gamma_k - \eta_k \big(R_k^{T}(\tilde{H}_{z}(z_{k}+R_{k}\gamma_k)) R_k\big)^{-1}R_k^{T}\nabla_{z} \tilde{f}(z_{k}+R_{k}\gamma_k)
			 \\
			 &=\gamma_k - \eta_k \big(R_k^{T}(H_{z}(z_{k}+R_{k}\gamma_k)+ \text{T}) R_k\big)^{-1}R_k^{T}\nabla_{z} \tilde{f}(z_{k}+R_{k}\gamma_k)
			,
\end{aligned}
\end{equation}
where the last two equations illustrate the computational complexity of the method, through the Hessian-vector product over the subspaces matrix $R_k$, and the low dimensional subspace Hessian inversion.
Here, the matrices  $H_{u}(v)$ and $\tilde{H}_{u}(v)$ denote the variable metric matrices reduced to $\nabla^{2}_{u}f(v)$ and $\nabla^{2}_{u}\tilde{f}(v)$ respectively in the Newton scheme. Also, $\text{T}$ denote the dampening matrix that ensures stability of the former saddle-point system, such that $\text{T}=\tau\left( \begin{array}{c|c}
    I & 0 \\
   \midrule
   0 & -I \\
\end{array}\right).$
Here $\eta_{k}$ is the step size commonly obtained via the line search procedure.
In the case of the Newton optimization in the subspace being computationally intensive (e.g. high dimensional subspace or prohibitive derivatives computation), a Quasi-Newton method can be deployed instead. In the saddle-point setting, Quasi-Newton alternatives that do \emph{not} enforce positive definiteness of the Hessian can be used, such as symmetric rank-one
(SR1) with usual handling of the update factors \cite{NocedWrite}.
We summarize the proposed sequential subspace optimization framework for saddle-point problems in Algorithm \ref{algo:sesop_saddle}.
\begin{center}
\begin{algorithm}
\caption{Sequential Subspace Saddle-point Optimization}\label{algo:sesop_saddle}
\begin{algorithmic}
\REQUIRE{$z_{0}=(x_0,y_0)$, $d$ the maximum subspace dimension, $K$ the maximum number of iterations, $\epsilon$ the machine precision}
\STATE {\bfseries Input:} $f:(\mathbb{R}^{M}\times\mathbb{R}^{N})\rightarrow \mathbb{R}$, initial point $(x_{0},y_{0})$
\STATE {\bfseries Output:} $(x_{final},y_{final})$
\STATE$\tau \in \mathbb{R}^{+}_{0},\nu\in(0,1)$
\STATE Initialize proximal centers $\bar{x}_{0}$, $\bar{y}_{0}$
\FOR[$\triangleright$ Outer loop of the method]{$k=0,1,...,K$}
\STATE \textbf{if} $\|\nabla{f}(x_{k},y_{k})\| < \epsilon$ \textbf{then}: \textbf{return} $(x_{k},y_{k})$
\STATE \textbf{if} $\|\nabla{\tilde{f}}(x_{k},y_{k})\| < \epsilon$ \textbf{then}: $\tau =\nu\tau$
\STATE Update $P_{k}$ and $Q_{k}$ with current gradients 
\STATE Set $t=0$ and $\gamma_{t}=0$
\WHILE[$\triangleright$ Inner loop of the method: Subspace optimization]{$\|\nabla_{\gamma}{f}(z_k+R_k\gamma_t)\| > \epsilon$}
\STATE $\bar{\gamma}=-\tilde{H}_{\gamma}^{-1}(z_k+R_k\gamma_t)\nabla_{\gamma}\tilde{f}(z_k+R_k\gamma_t)$, eq. (\ref{eq:sesop_minimax_tau}).
\STATE Find inner step size $\eta_{in}$ following Algorithm \ref{algo:line_search} over subspace objective
\STATE Set $\gamma_{t+1}=\gamma_{t}+\eta_{in}\bar{\gamma}$ and $t=t+1$
\ENDWHILE
\STATE Find outer step size $\eta_{out}$ following Algorithm \ref{algo:line_search}
\STATE Update $z_{k+1}=z_{k}+\eta_{out}R_{k}\gamma_{t}$
\STATE \textbf{if} $dim(P_k) > d-1$ \textbf{then}: Remove the oldest direction from $P_k$ and $Q_k$
\STATE Update $P_{k+1}$ and $Q_{k+1}$ with search steps and/or gradients
\STATE Update proximal centers $\bar{x}_{k}$, $\bar{y}_{k}$ (e.g with $x_{k}$ and $y_{k}$)
\ENDFOR
\STATE \textbf{return} $(x_{K},y_{K})$
\end{algorithmic}
\end{algorithm}
\end{center}
\section{Experimental Results}
\label{others}
To assess the performance of the proposed method, and to demonstrate its efficacy, we performed several experiments.
GDA refers to the gradient method \cite{rockafellar1976monotone}, DAS refers to the dual averaging scheme of \cite{nesterov2009primaldual}, OGDA is the optimistic gradient method \cite{daskalakis2017training}, and EGDA refers to the extrinsic gradient method \cite{yadav2017stabilizing}.
All the methods but DAS are implemented using the proposed backtracking line search for improved convergence. When line search did not converge for EGDA, we searched for optimal step size. 
In the following, the proximal centers of the proposed method are set to previous point (i.e. $x_{k-1}$ and $y_{k-1}$).
In all the presented experiments, the subspace is populated by the current gradient ($m=n=1$), by previous gradient ($m=n=2$), and by previous search directions ($m=n\geq 3$).
The machine precision $\epsilon$ is set to single precision $10^{-8}$ and the maximum number of inner optimization iterations is limited to 10.
In all the figures $k$ denote the iteration number. Unless stated otherwise all the methods use the same oracle.
\subsection{Quadratic saddle-point Problem}
We consider the following quadratic saddle-point problem
\begin{equation} \label{eq:saddle_quadratic}
\begin{aligned}
\min_x \max_y \ \frac{1}{2}(x^T A_{x} x + y^T A_{y} y) + x^{T}Cy + b_{x}^{T}x+b_{y}^{T}y, 
\end{aligned}
\end{equation}
where the matrices $A_x, A_y, C$ are generated from the normal distribution and have pre-defined condition numbers.
Namely, we generate a standard Gaussian matrix with i.i.d. entries, perform its SVD and substitute diagonal singular values with an array of log-uniform random values in predefined range.
The dimension of the optimization problem is set to $M=1500,N=500$.
We plot the distance to optimum (leftmost plot) and the norm of the gradient (second plot from left).
Also, we show the impact of the condition number of the block-matrices in $A$ on the mean convergence rate $K^{-1}\sum_{k}^{K}{\|z_{k+1}-z^{*}\|}/{\|z_{k}-z^{*}\|}$ (rightmost plot).
Here $x$-axis represents the inverse condition number $\kappa^{-1}$.
These first three plots are obtained with a \emph{three} dimensional subspace for each variable.
Finally, we show the effect of the subspace dimension on the convergence of the proposed method (third from left). 
For this experiment we add time comparison with the Chambolle-Pock algorithm (CP) with step size set according to the interaction term matrix to ensure proven convergence \cite{chambolle2011first}.
For fairness, we do not assume any closed form solution is given, all the derivatives are computed at each iteration using the same automatic differentiation tool for all the methods.

\begin{figure}[h]
\begin{center}
\subfigure{\includegraphics[trim={0 70 0 0}, width=0.24\columnwidth]{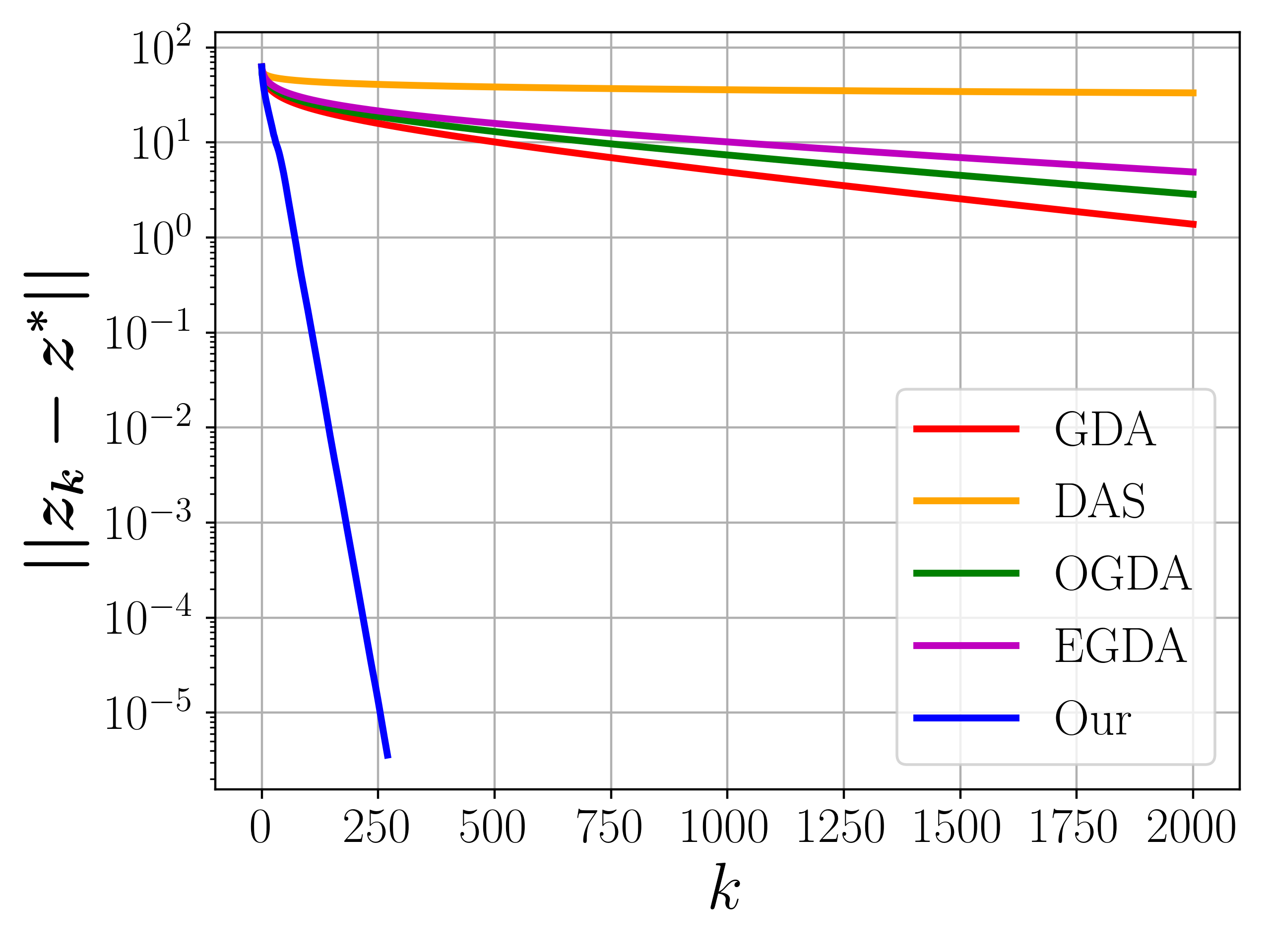}}
\subfigure{\includegraphics[trim={0 70 0 0}, width=0.24\columnwidth]{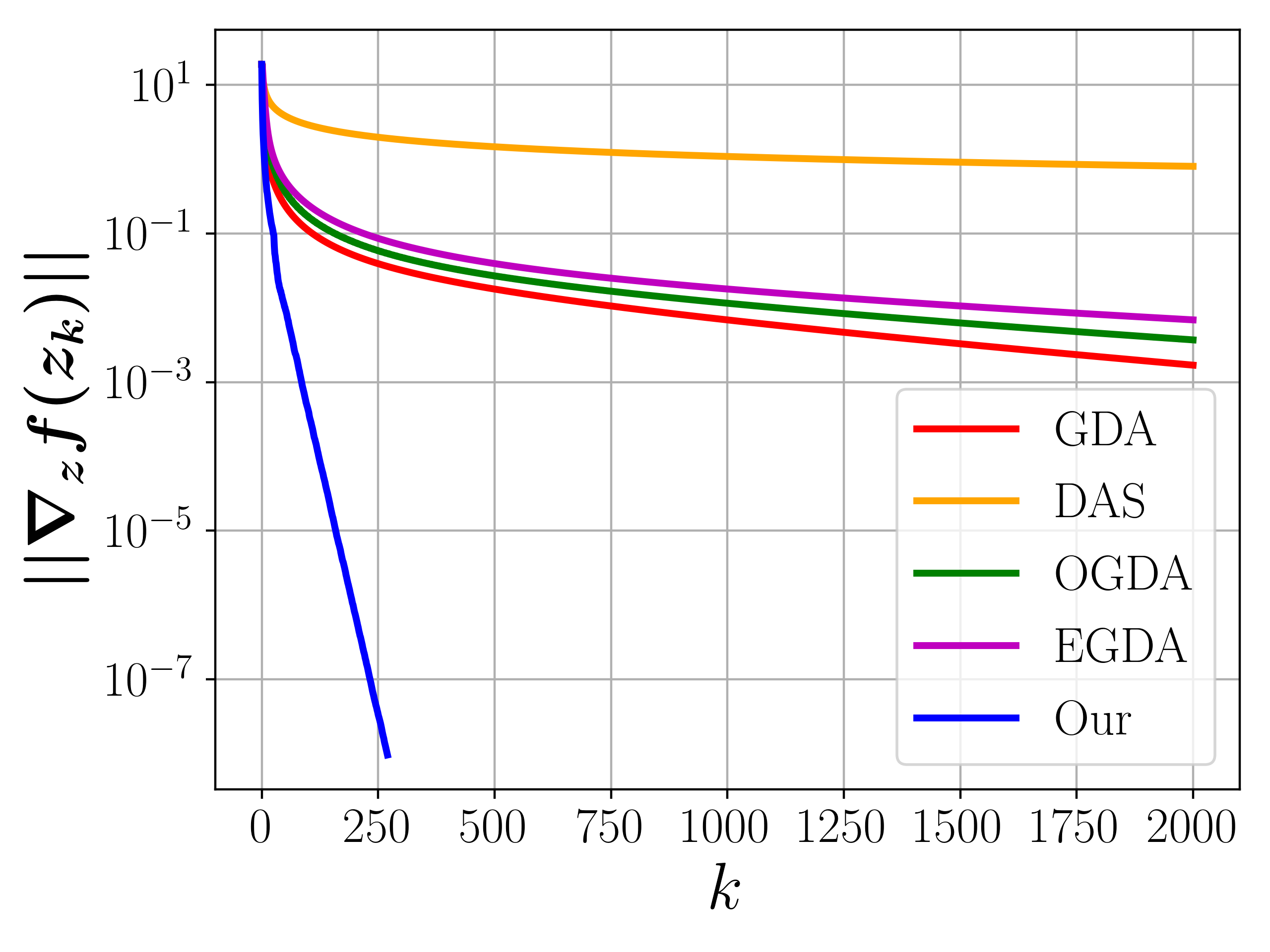}}
\subfigure{\includegraphics[trim={0 70 0 0}, width=0.24\columnwidth]{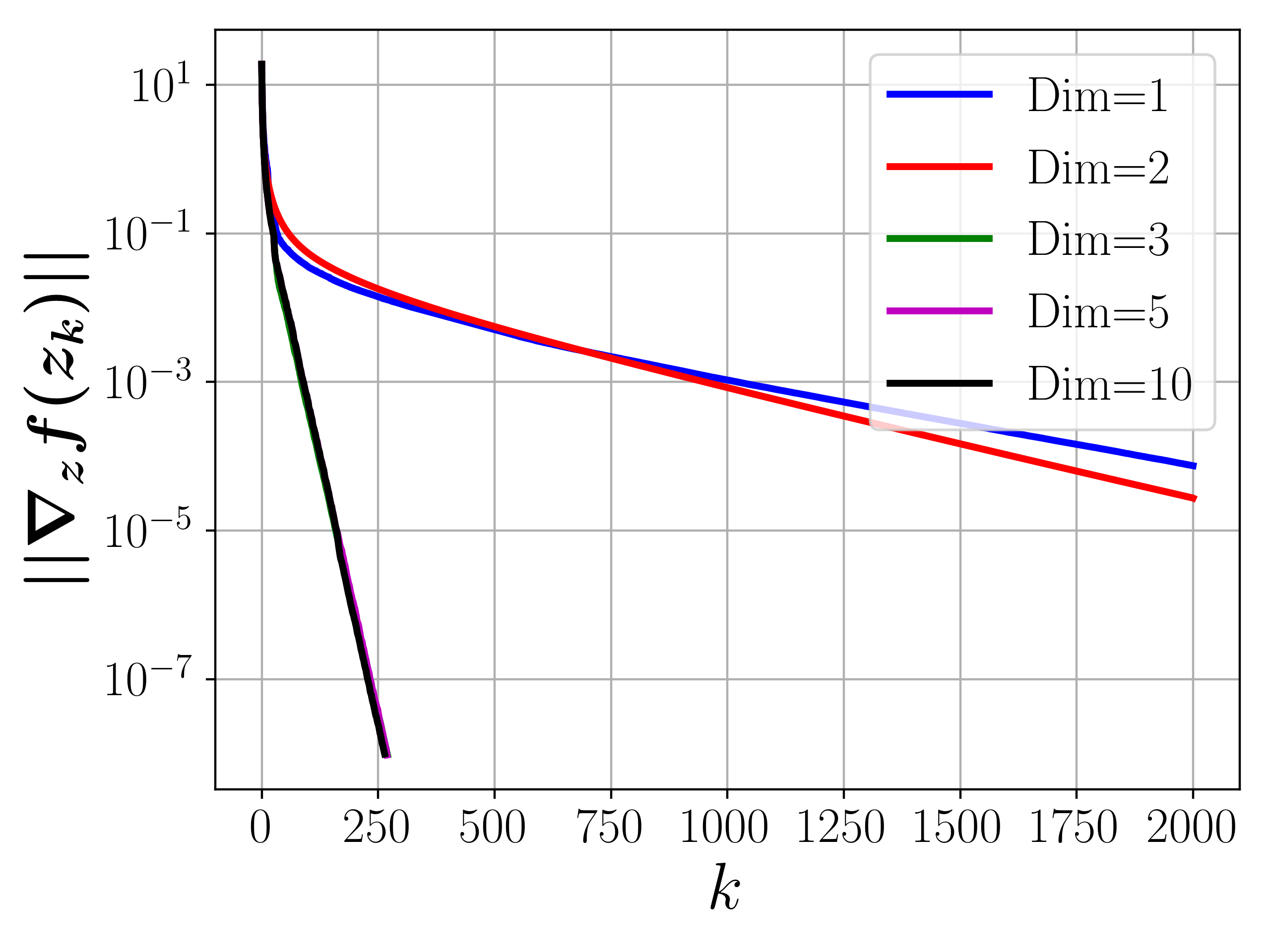}}
\subfigure{\includegraphics[trim={0 70 0 0}, width=0.24\columnwidth]{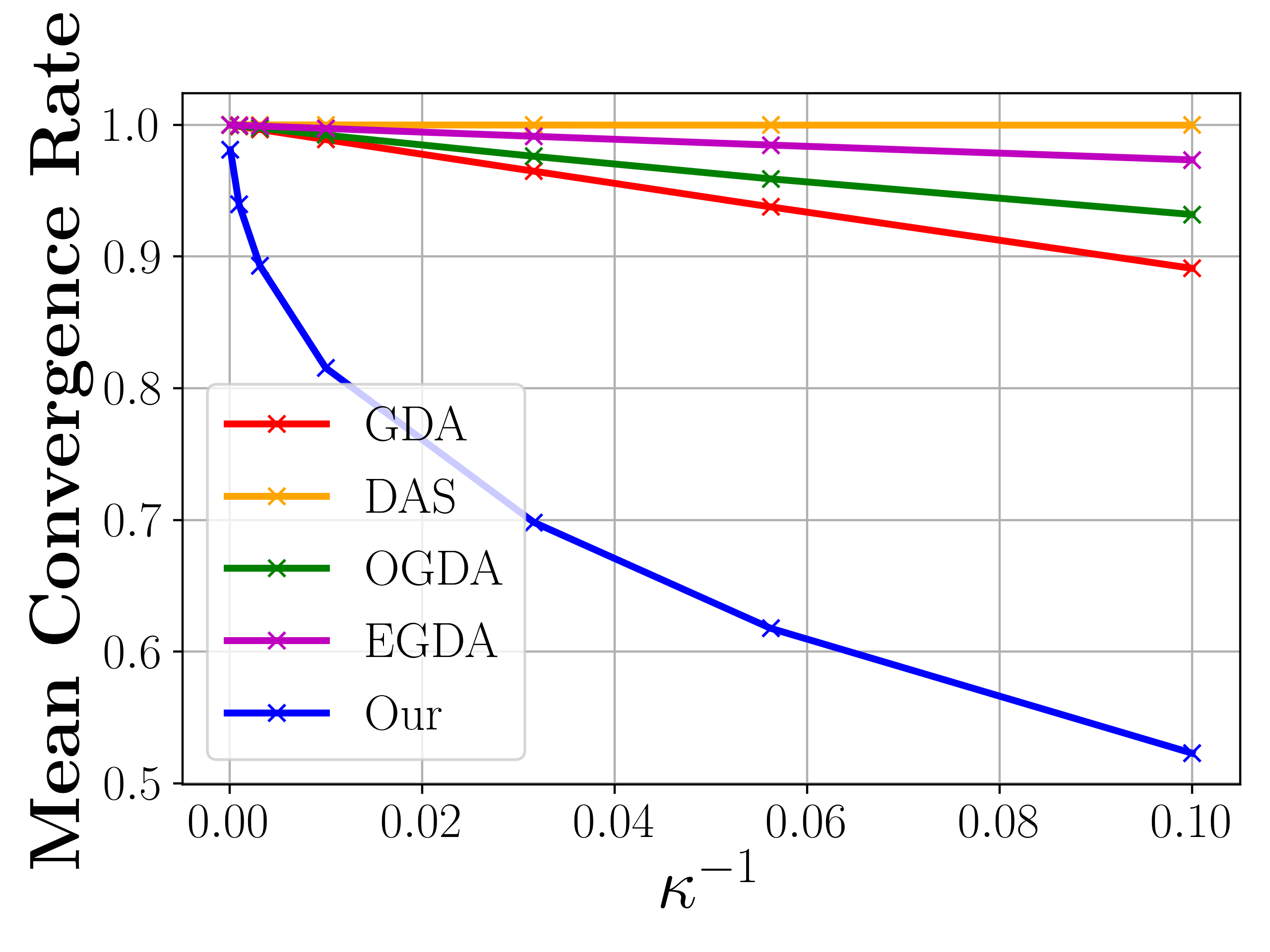}}
\caption{Separable quadratic saddle point problem}
\label{fig:separable}
\end{center}
\end{figure}

In the first experiment presented in Figure \ref{fig:separable}, we consider a separable problem with $A_x \succ 0$, $A_y\prec 0$ and $C=0$. The two matrices are conditioned with the condition numbers $\kappa(A_x)=10^3,\kappa(A_y)=10^2$.
Figure \ref{fig:separable} shows that the proposed approach keeps its manifold expansion property throughout the last search direction ($m=n\geq3$), and therefore converges extremely fast to the solution. In contrast, the convergence of other methods is very slow. This is due to the 
difficulty of the gradient method to converge in ill-conditioned scenarios, and due to the unified step size for the primal and the dual directions.

\begin{figure}[h]
\begin{center}
\setcounter{subfigure}{0}
  \subfigure{\includegraphics[trim={0 70 0 0}, width=0.24\columnwidth]{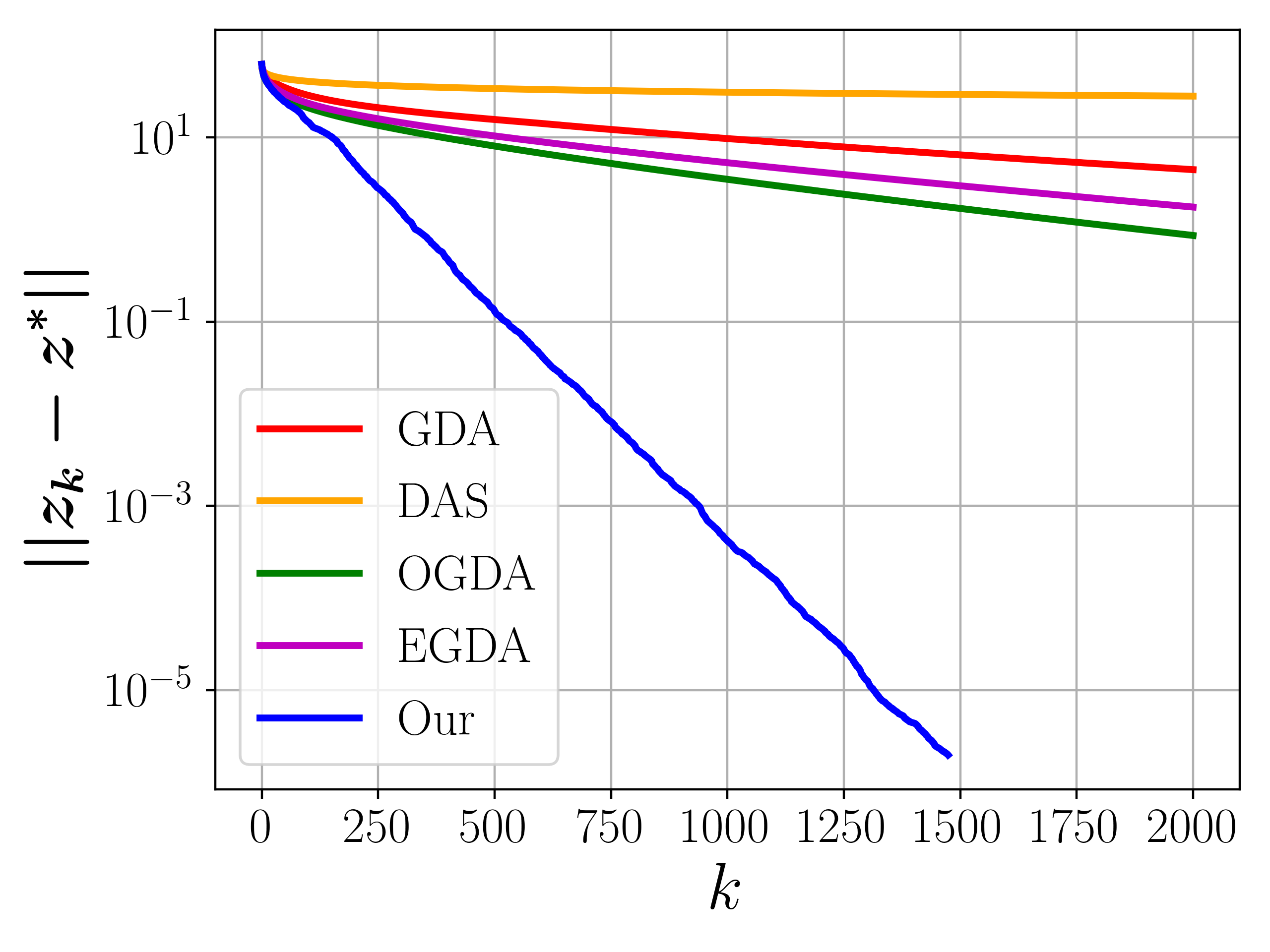}}
  \subfigure{\includegraphics[trim={0 70 0 0}, width=0.24\columnwidth]{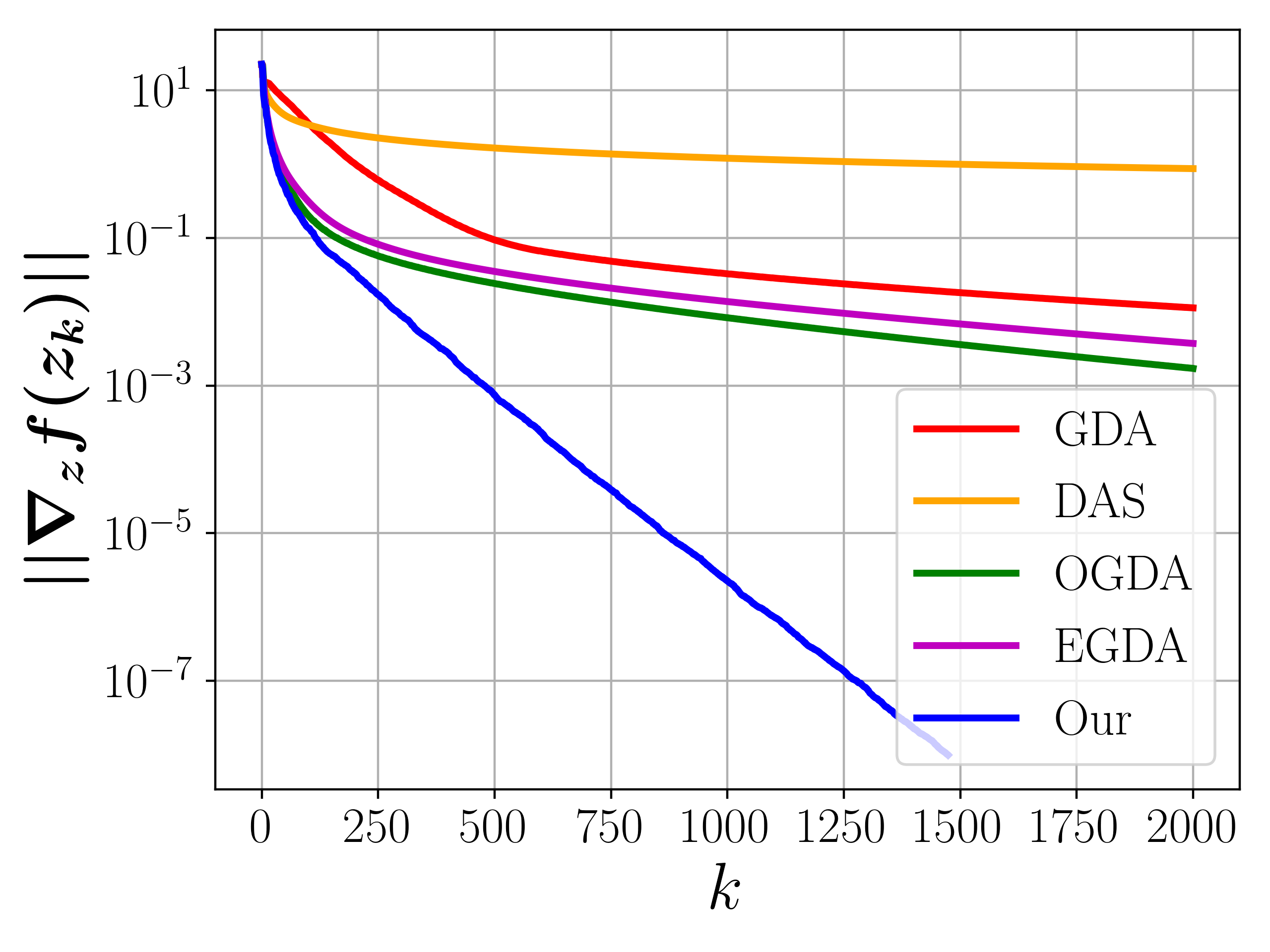}}
  \subfigure{\includegraphics[trim={0 70 0 0}, width=0.24\columnwidth]{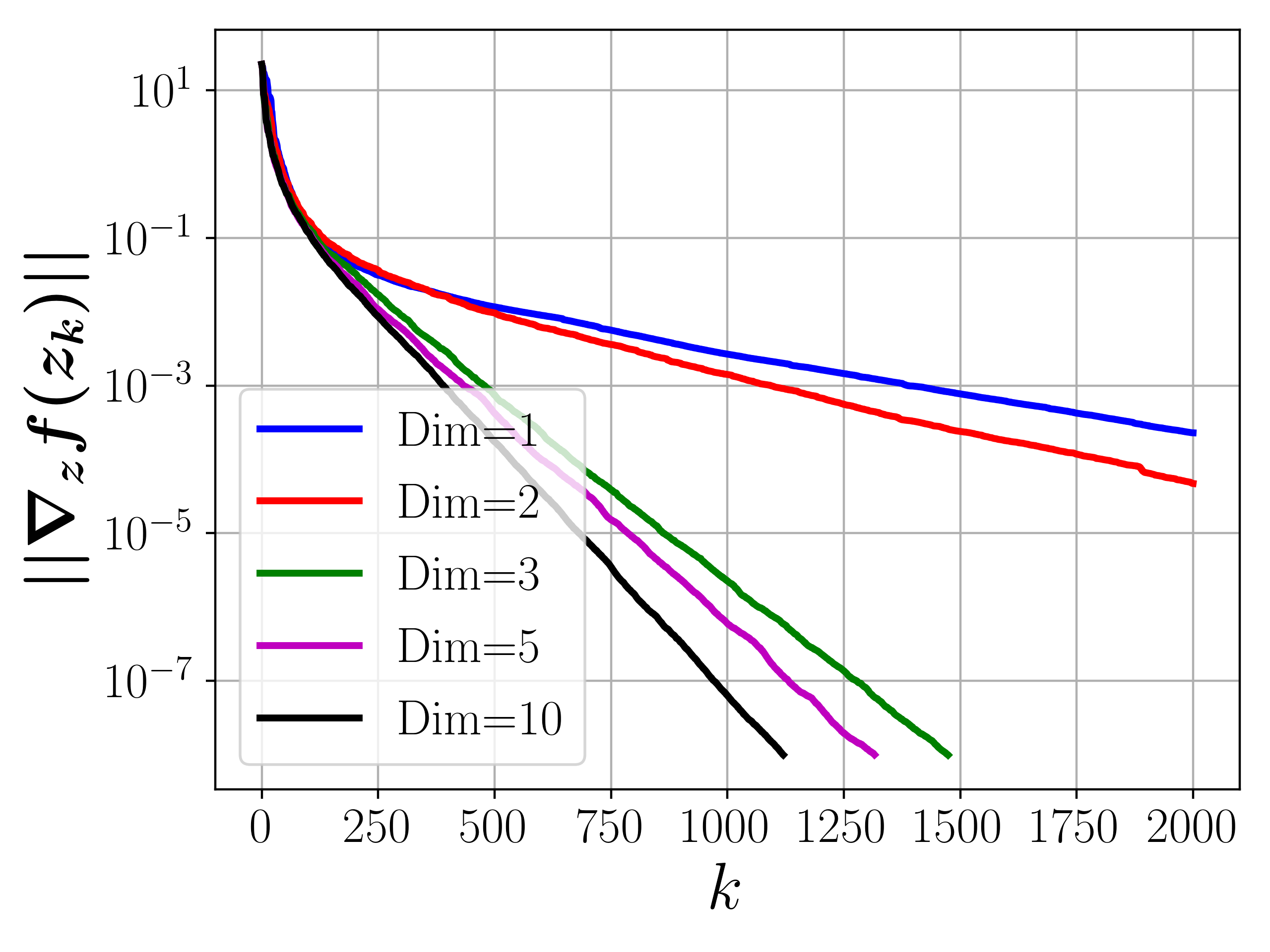}}
  \subfigure{\includegraphics[trim={0 70 0 0}, width=0.24\columnwidth]{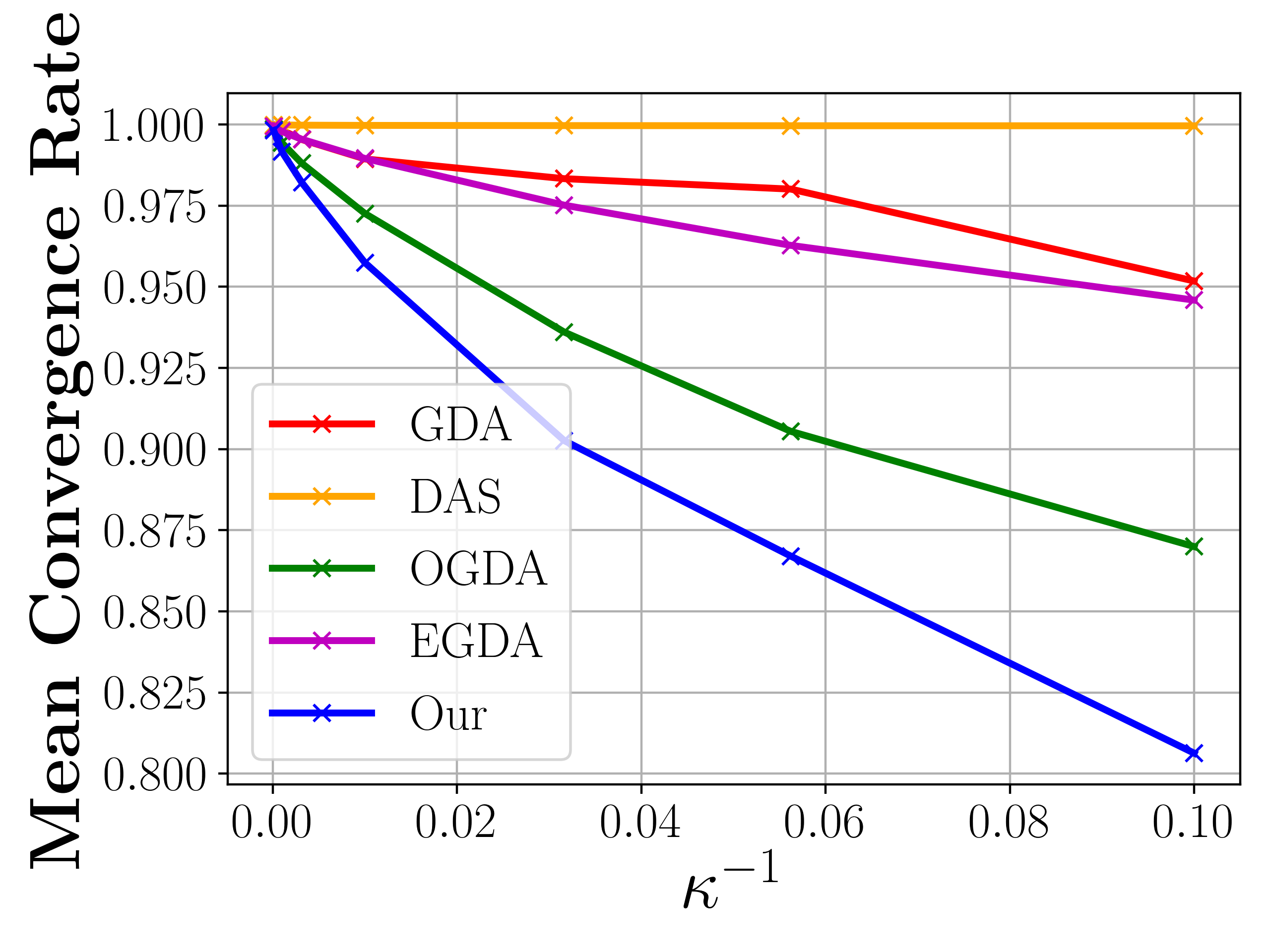}}
\caption{Stable quadratic saddle point problem}
\label{fig:stable}
\end{center}
\end{figure}

\begin{figure}[h]
\begin{center}
\setcounter{subfigure}{0}
  \subfigure{\includegraphics[trim={0 70 0 0}, width=0.24\columnwidth]{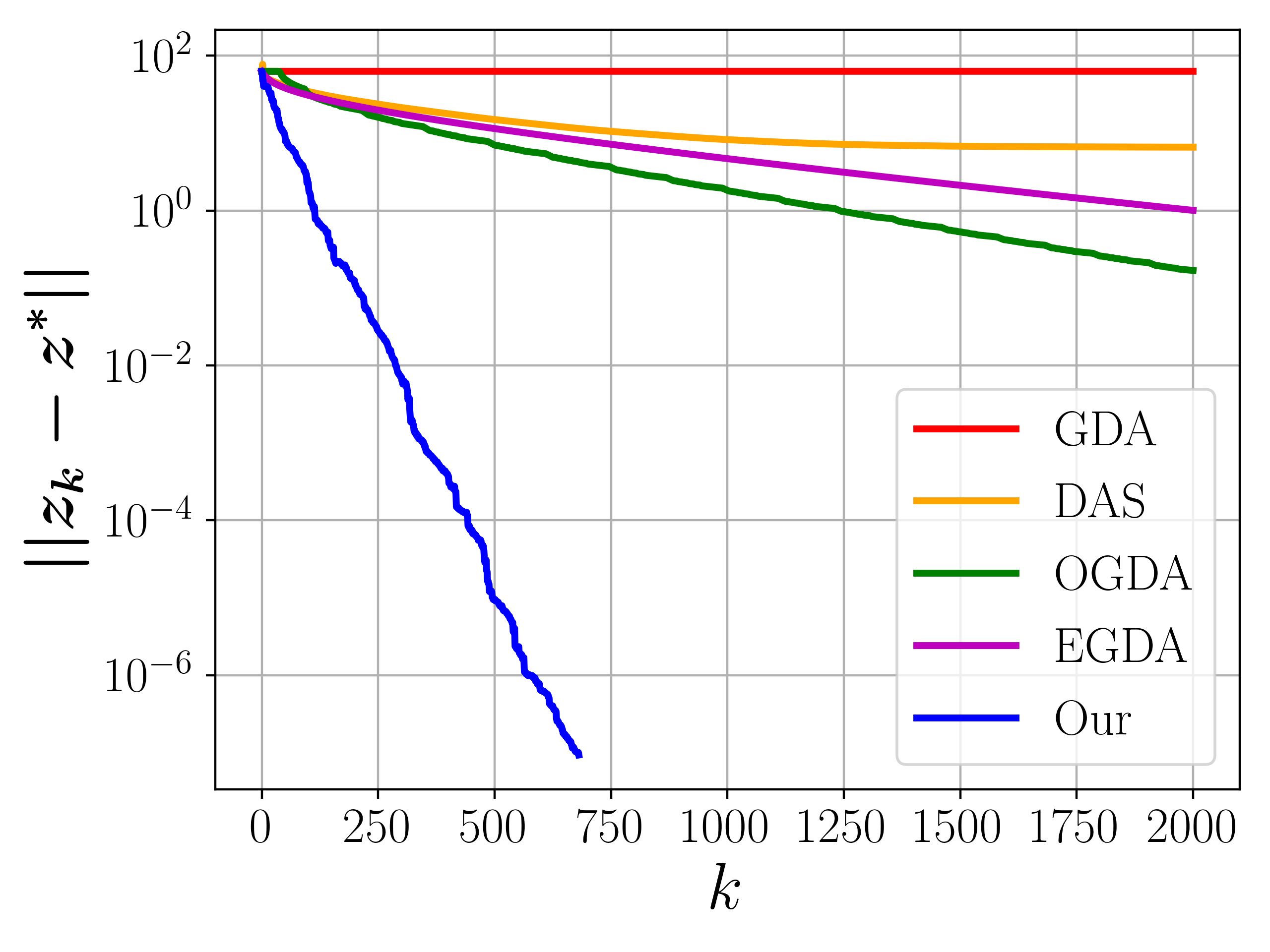}}
  \subfigure{\includegraphics[trim={0 70 0 0}, width=0.24\columnwidth]{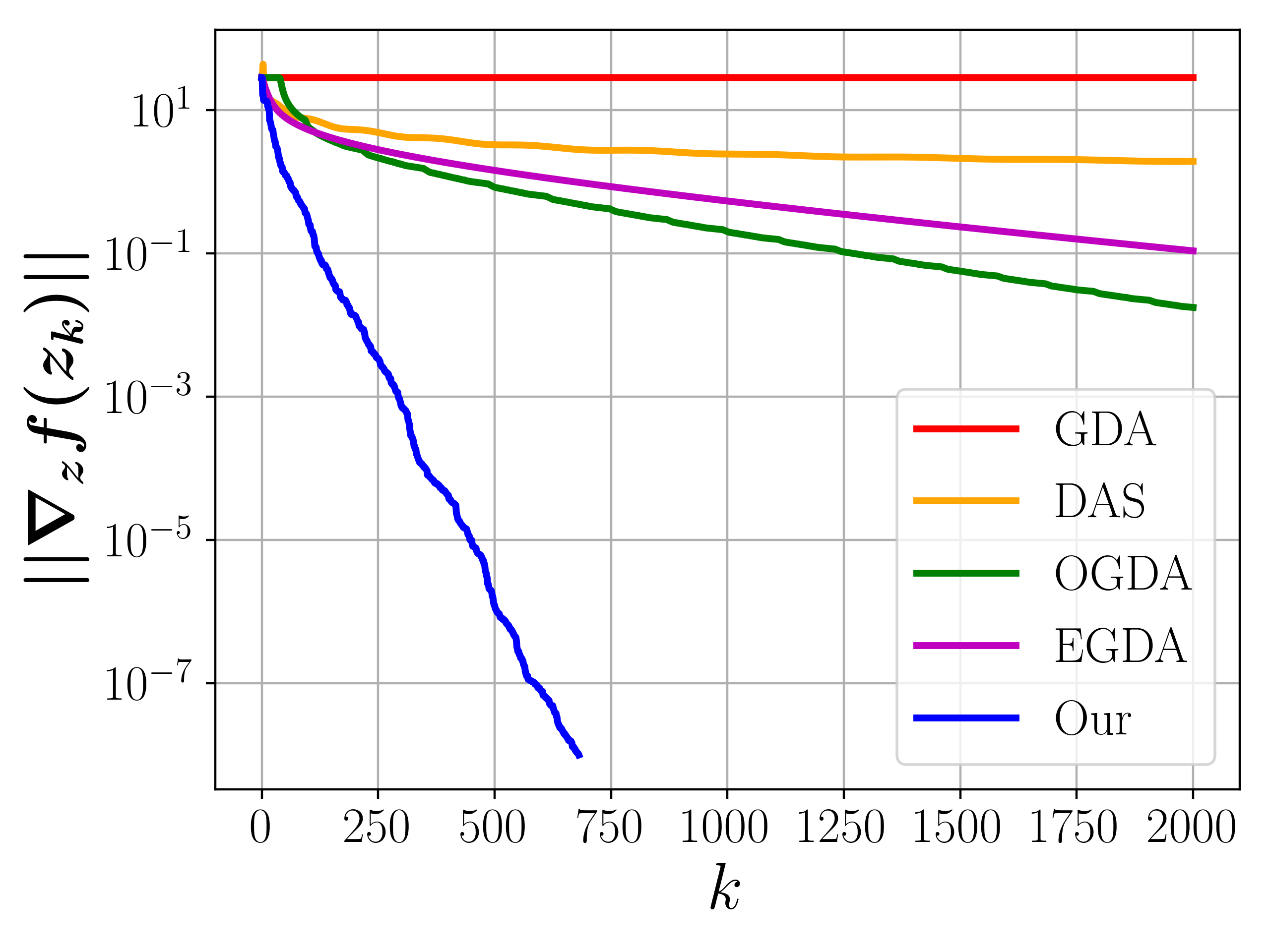}}
  \subfigure{\includegraphics[trim={0 70 0 0}, width=0.24\columnwidth]{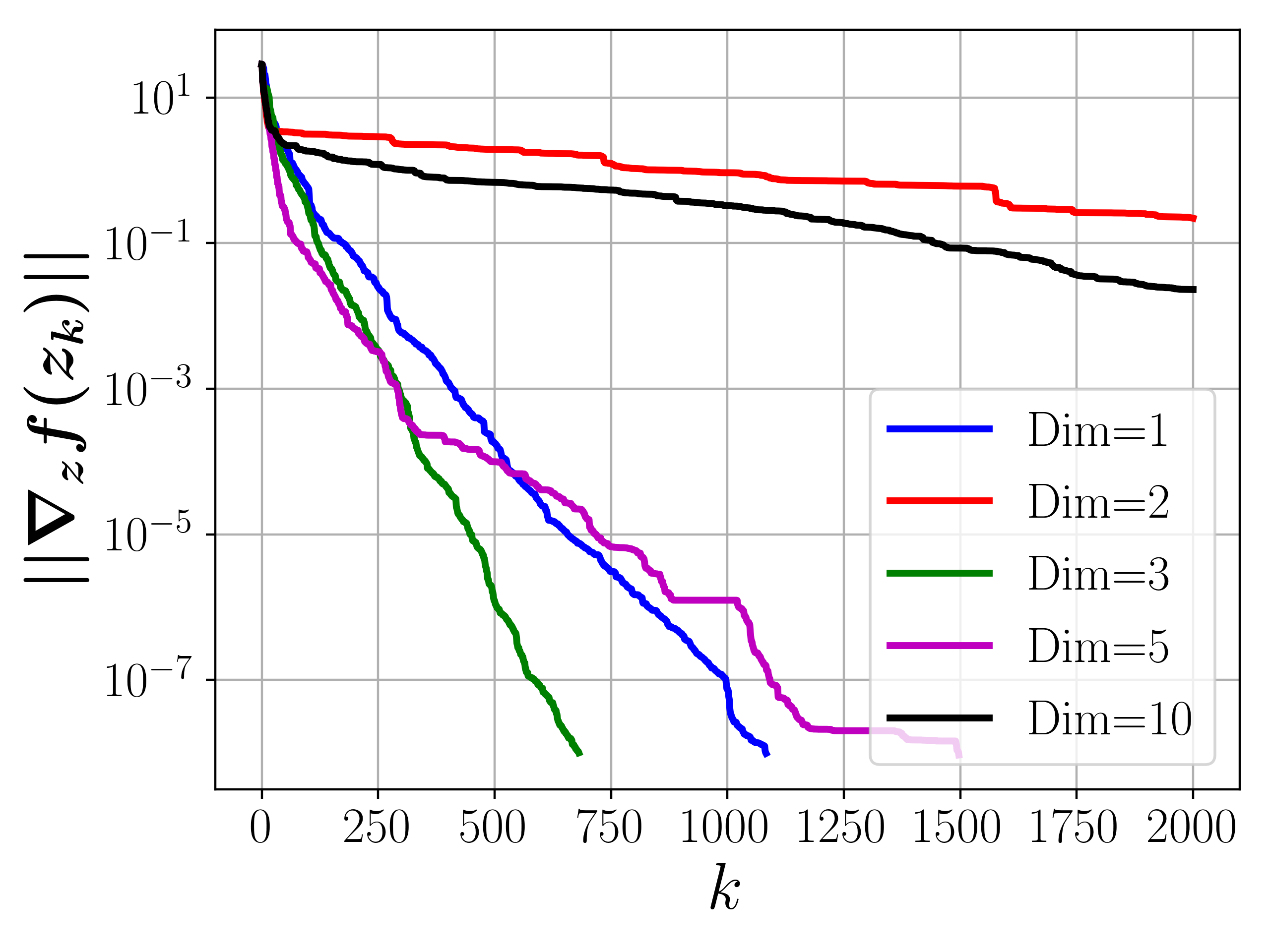}}
  \subfigure{\includegraphics[trim={0 70 0 0}, width=0.24\columnwidth]{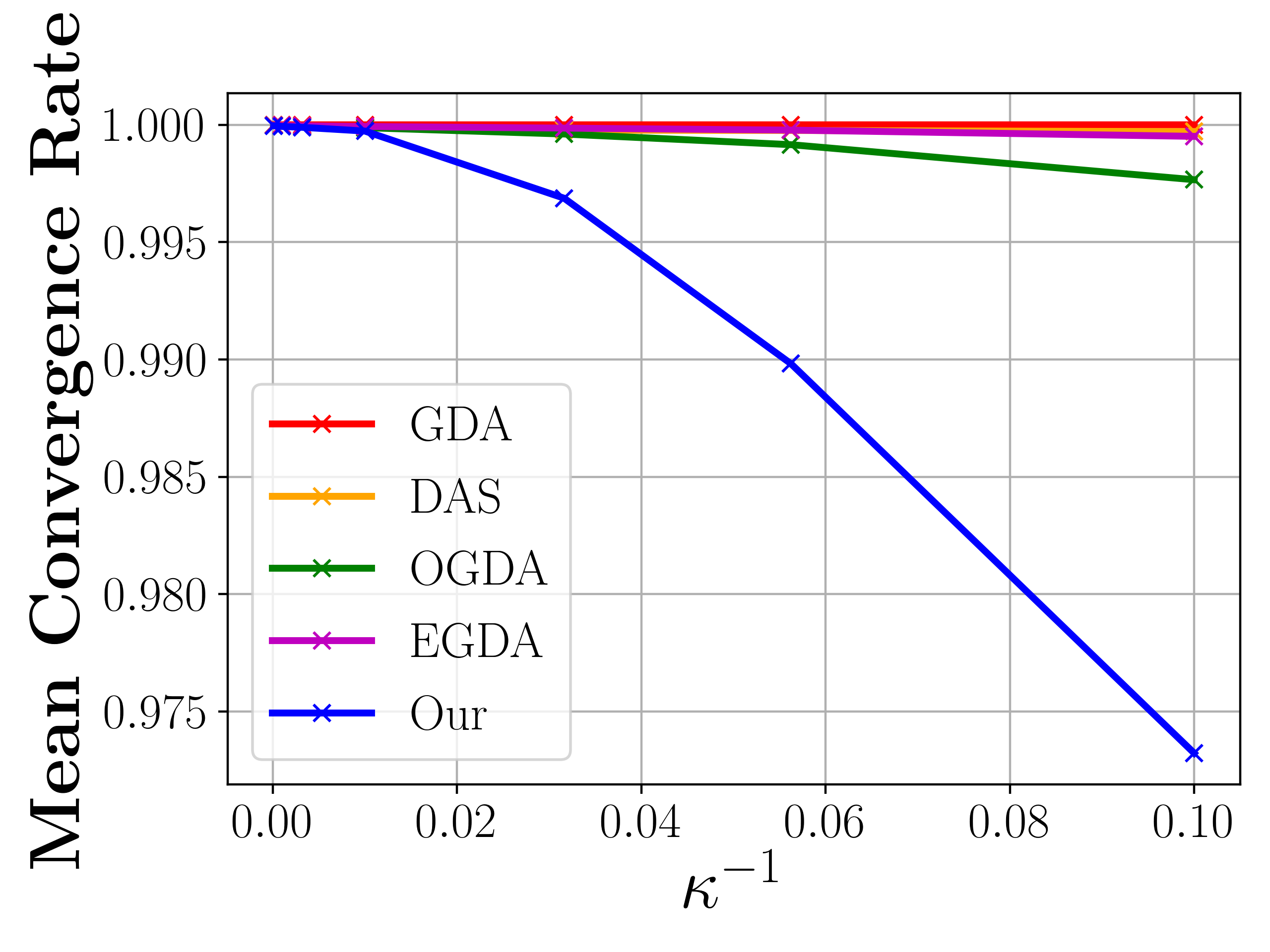}}
\caption{Unstable quadratic saddle point problem}
\label{fig:bilinear}
\end{center}
\end{figure}

In the second experiment presented in Figure \ref{fig:stable}, we consider the stable quadratic saddle point problem with $A_x \succ 0$, $A_y\prec 0$ and $C$ to be full rank matrices. Here, all the  block-matrices are conditioned with condition number $\kappa(A_x)=10^3,\kappa(A_y)=10^2,\kappa(C)=10^3$.
We observe the superiority of the proposed method, while the advantage of using more directions is clear in handling the interaction matrix $C$, as compared to gradient based methods.
In the last experiment, we deploy the bi-linear game problem, where $A_x=0$, $A_y= 0$ and $C$ is a full rank matrix, such that $\kappa(C)=10^2$. Here $M=N=1000$ so there exist only one solution to the equivalent system of linear equations. The results of the experiments are presented in Figure \ref{fig:bilinear}.
The proposed method performs significantly better, as compared to other first-order approaches.
In the bi-linear case, we can see that increasing the size of the subspace is not necessarily beneficial.  
We show the superiority of the method in term of computational time for the different settings in Table \ref{table:time_runs_quadratic}.
\begin{table}\label{table:time_runs_quadratic}
\begin{center}
    \begin{tabular}{| l | l | l | l | l | l | l |}
    \hline
{Setting} &   {GDA}     &   {DAS}     &   {OGDA}        &   {EGDA}    &   {CP}      &      \textbf{ Our}     \\ \hline
{Separable}   &   42.5   &   -     &   73.1    &   60.7     &   15.9     &   \textbf{5.9}     \\ \hline
{Stable}      &   144.1     &   -     &   49.3     &   55.4     &   -     &   \textbf{38.1}     \\ \hline
{Unstable}    &   $\infty$     &   -     &   34.3     &    33.6    &   -     &   \textbf{14.7}     \\ \hline
    \end{tabular}
\end{center}
        \caption{Mean computation time in seconds of the presented methods until convergence threshold is reached, for the different quadratic settings. $\infty$ denotes non-convergence and '-' denotes slower convergence than our method by at least factor 30.}
    \vspace{-0.5cm}
\end{table}
\subsection{Constrained Optimization: ADMM}
\begin{figure}[h]
\begin{center}
\setcounter{subfigure}{0}
  \subfigure{\includegraphics[trim={0 70 0 0}, width=0.24\columnwidth]{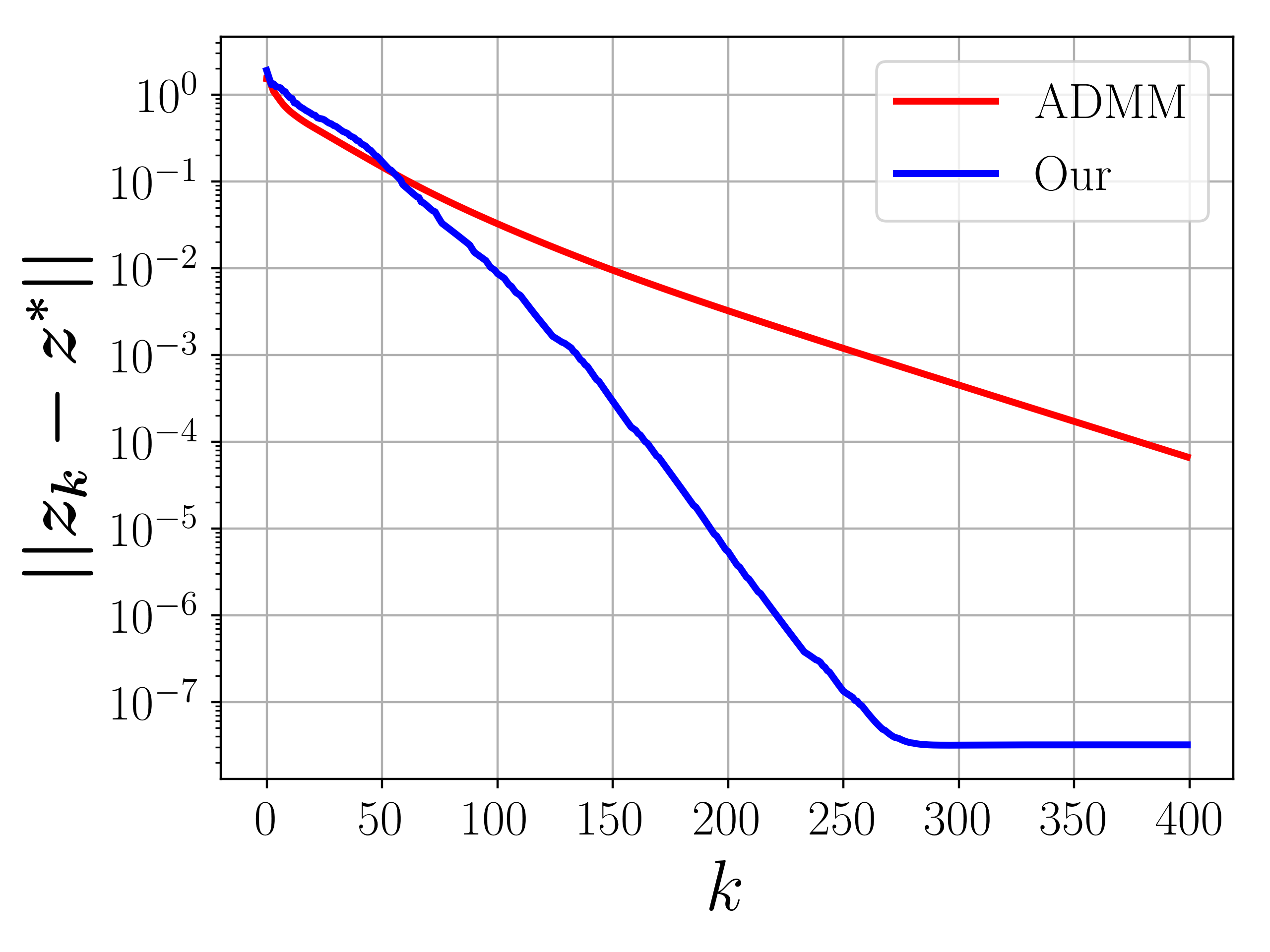}}
  \subfigure{\includegraphics[trim={0 70 0 0}, width=0.24\columnwidth]{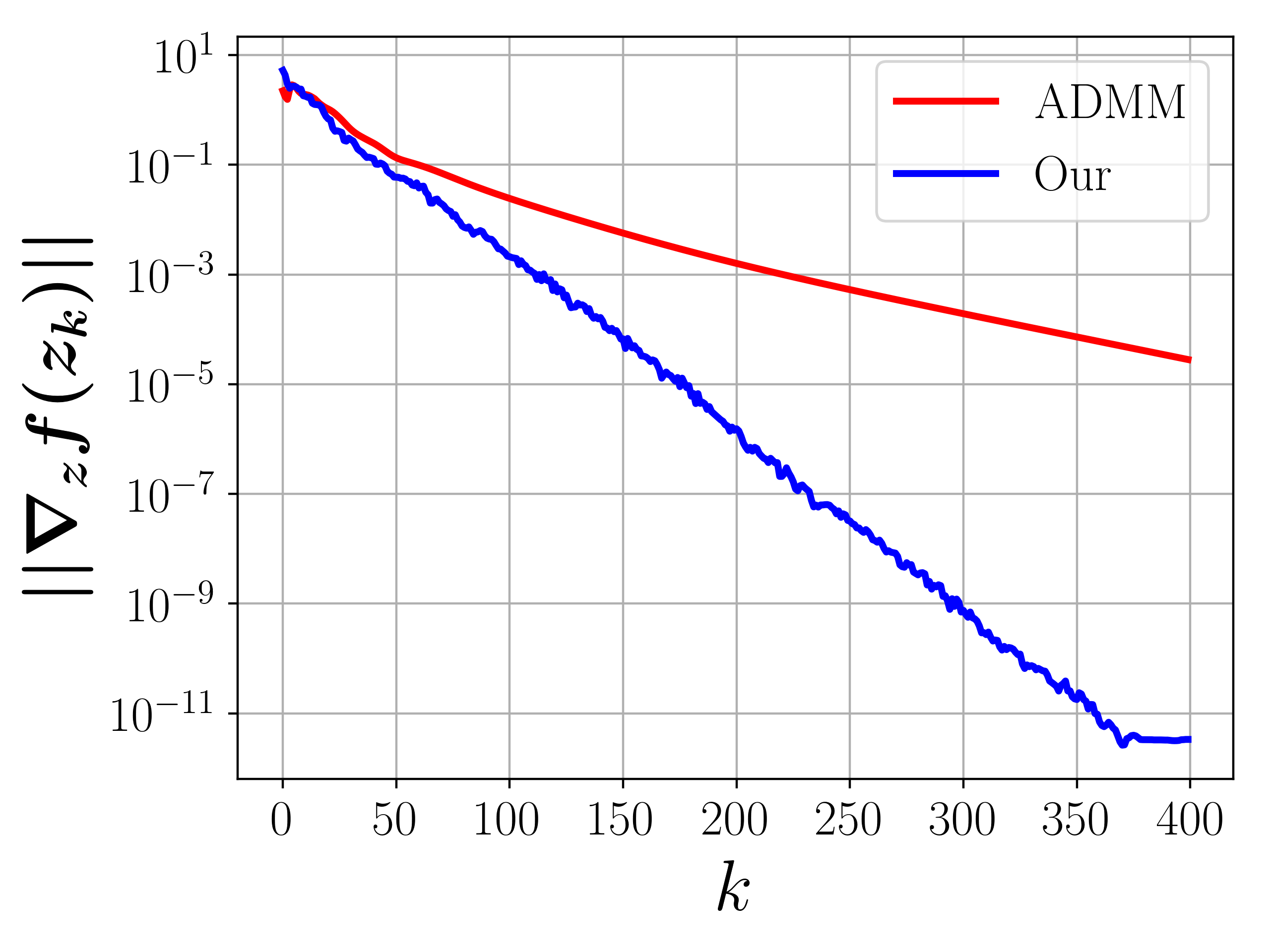}}
\caption{Subspace Optimization boosting via ADMM Directions }
\label{fig:admm}
\end{center}
\end{figure}

The alternating direction method of multipliers (ADMM) is a powerful algorithm that solves convex optimization problems by splitting them into simpler problems that are easier to handle.
In this experiment, we consider the smooth Lasso regression problem, where the smoothness is enforced in order to allow efficient continuous optimization. 
The problem is defined as follows
\begin{equation} \label{eq:saddle_augmented}
\begin{aligned}
& \underset{x}{\min}
& &  \frac{1}{2}\|Ax-b\|^{2}+\lambda \sum_{j}\varphi_{s}(x_{j}),
\end{aligned}
\end{equation}
with $ \lambda \in \mathbb{R^{+}}$. Here $x_j$ denotes the $j^{th}$ component of vector $x$, and $\varphi_{s}(t)$ denotes the scalar non-linearity that implements the smooth convex approximation $\sum_{j}\varphi_{s}(x_{j})$ of the $\ell_{1}$ norm such that 
\begin{equation} \label{eq:approx_l1}
\varphi_{s}(t)=|t| - s \text{ln}(1+|t|/s), s\in (0,\infty),
\end{equation}
where the scaling factor $s$ defines the degree of smoothness.
This choice of $\varphi_{s}(t)$ yields well defined shrinkage \cite{elad2007coordinate}.
The original ADMM algorithm can be summarized in the following three steps: minimization in former primal variable $x$, minimization in separable variable $w$, and update of the dual variable $y$ \cite{Boyd_alternating_2011}.
We can reformulate the Lasso setting as a saddle-point problem of the augmented Lagrangian
\begin{equation} \label{eq:saddle_admm}
\begin{aligned}
\underset{x,w}{\min}\ \underset{y}{\max} \ \Big\{
\frac{1}{2}\|Ax-b\|^2+ \lambda\sum_{j}\varphi_{s}(w_{j}) + y^{T}(x-w) + \frac{\rho}{2}\|x-w\|^2\Big\}
\end{aligned}
\end{equation}
where $\rho$ denotes the penalty parameter.
In Figure \ref{fig:admm} we present the convergence results of the ADMM method with smoothing constant $s=10^{-3}$, versus the proposed subspace method boosted by the ADMM directions populating the subspace matrices.
The data setting is the same as in \cite{Boyd_alternating_2011}, Section (11.1).
The boosting obtained by the proposed approach is significant both in speed and accuracy.
\subsection{Generative Adversarial Networks}
\begin{figure}[h]
\begin{center}
    \setcounter{subfigure}{0}
    \subfigure{\includegraphics[trim={0 70 0 0}, width=0.24\columnwidth]{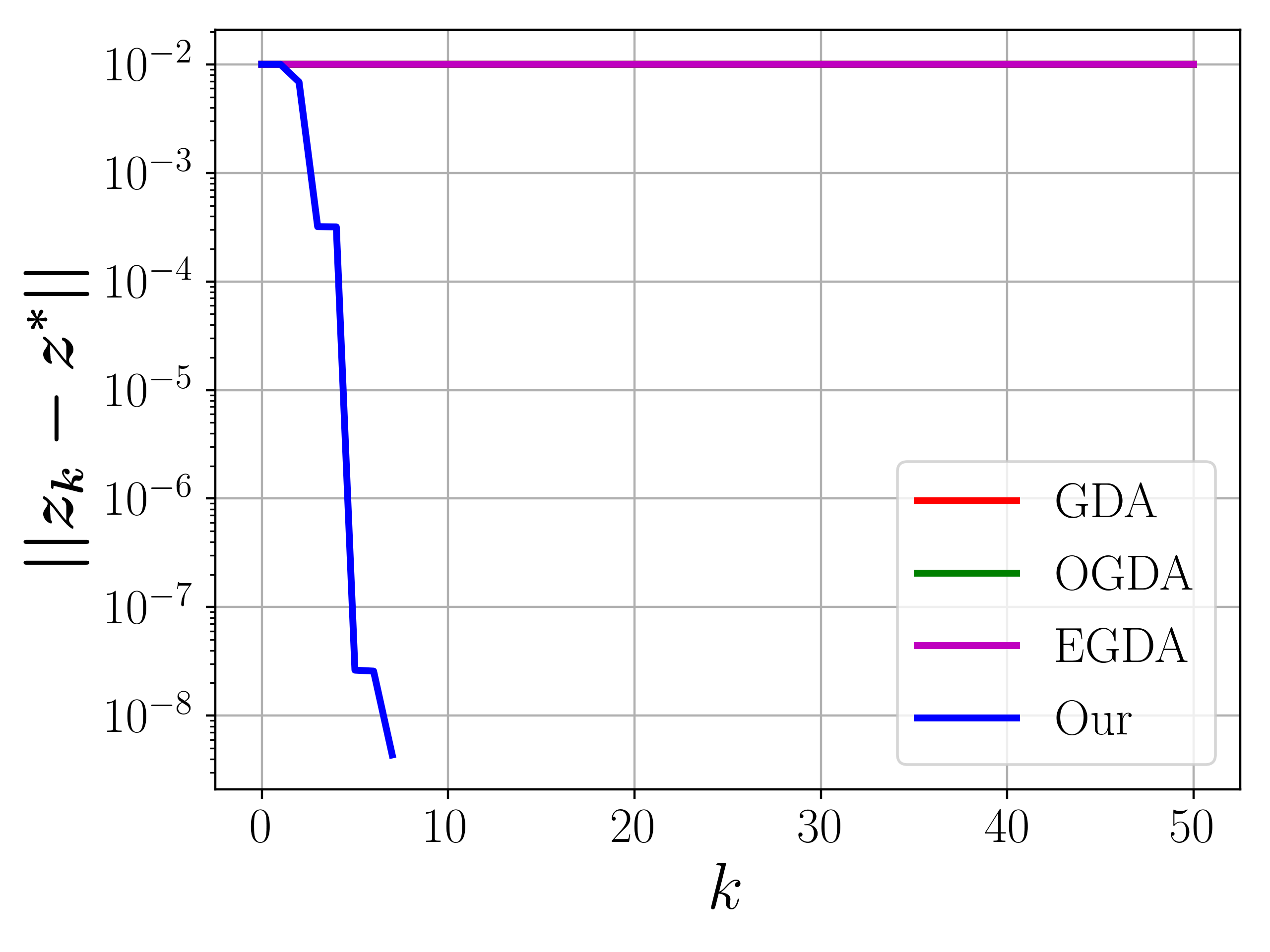}}
    \subfigure{\includegraphics[trim={0 70 0 0}, width=0.24\columnwidth]{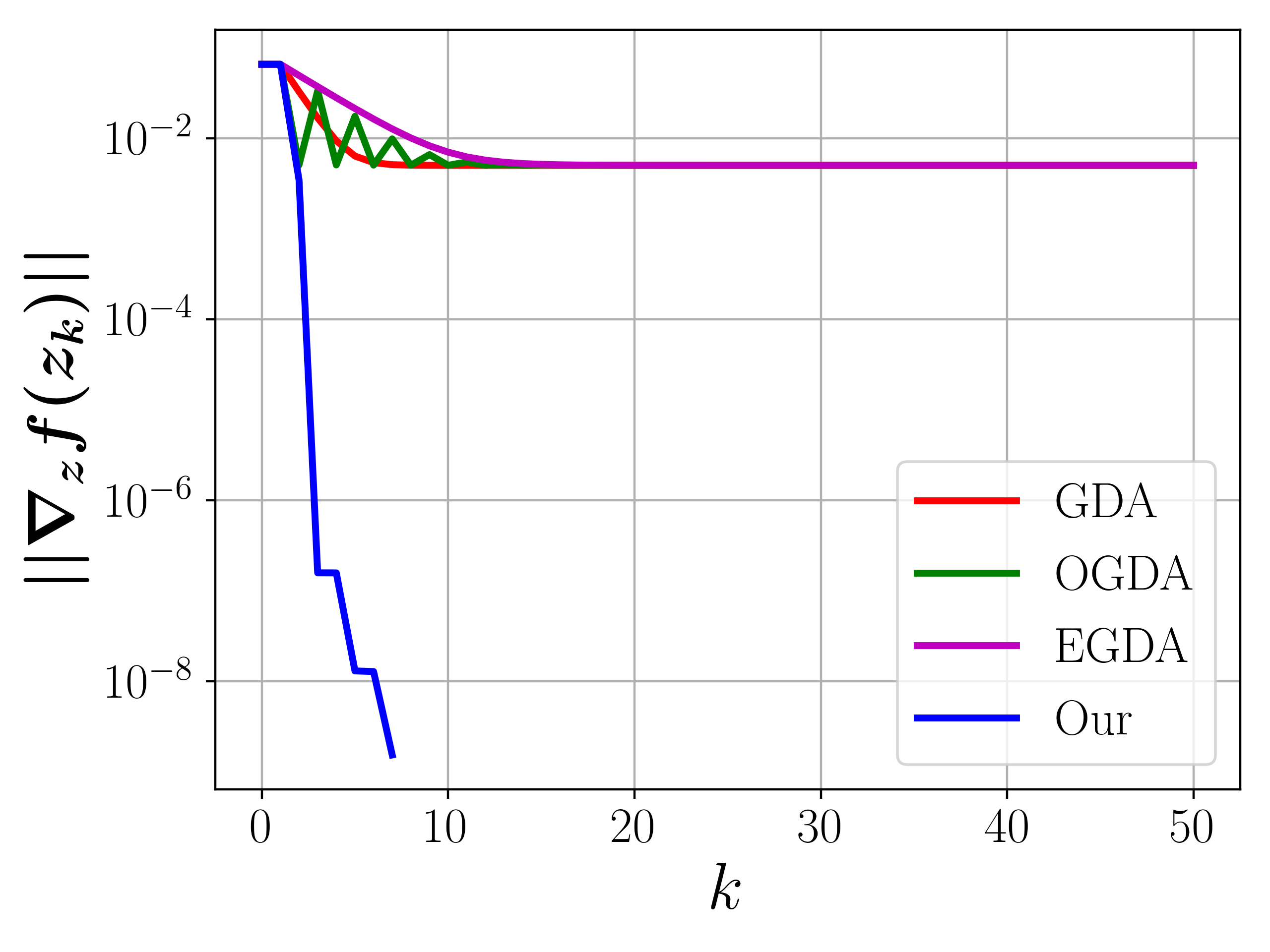}}
    \subfigure{\includegraphics[trim={0 70 0 0}, width=0.24\columnwidth]{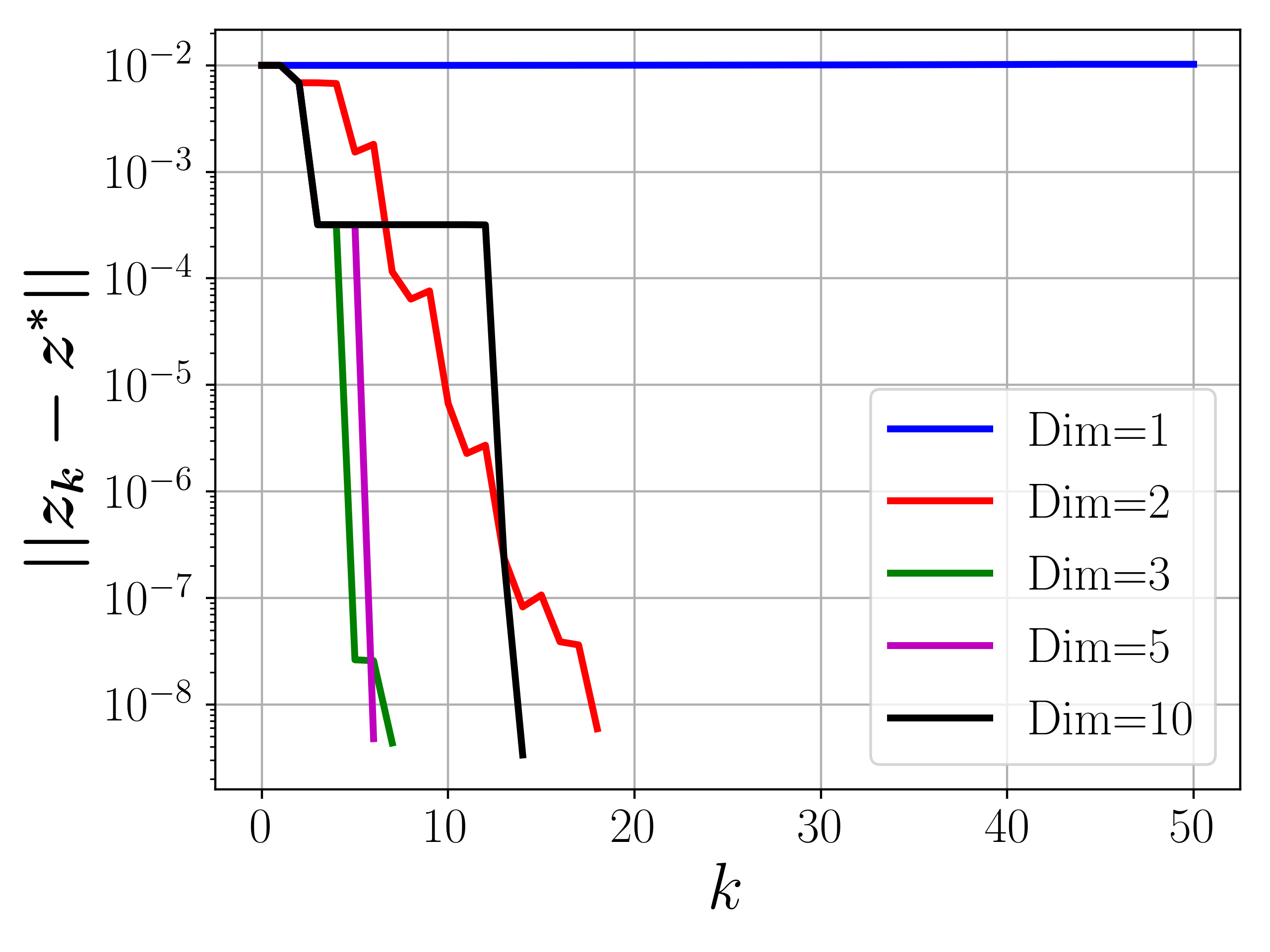}}
    \subfigure{\includegraphics[trim={0 70 0 0}, width=0.24\columnwidth]{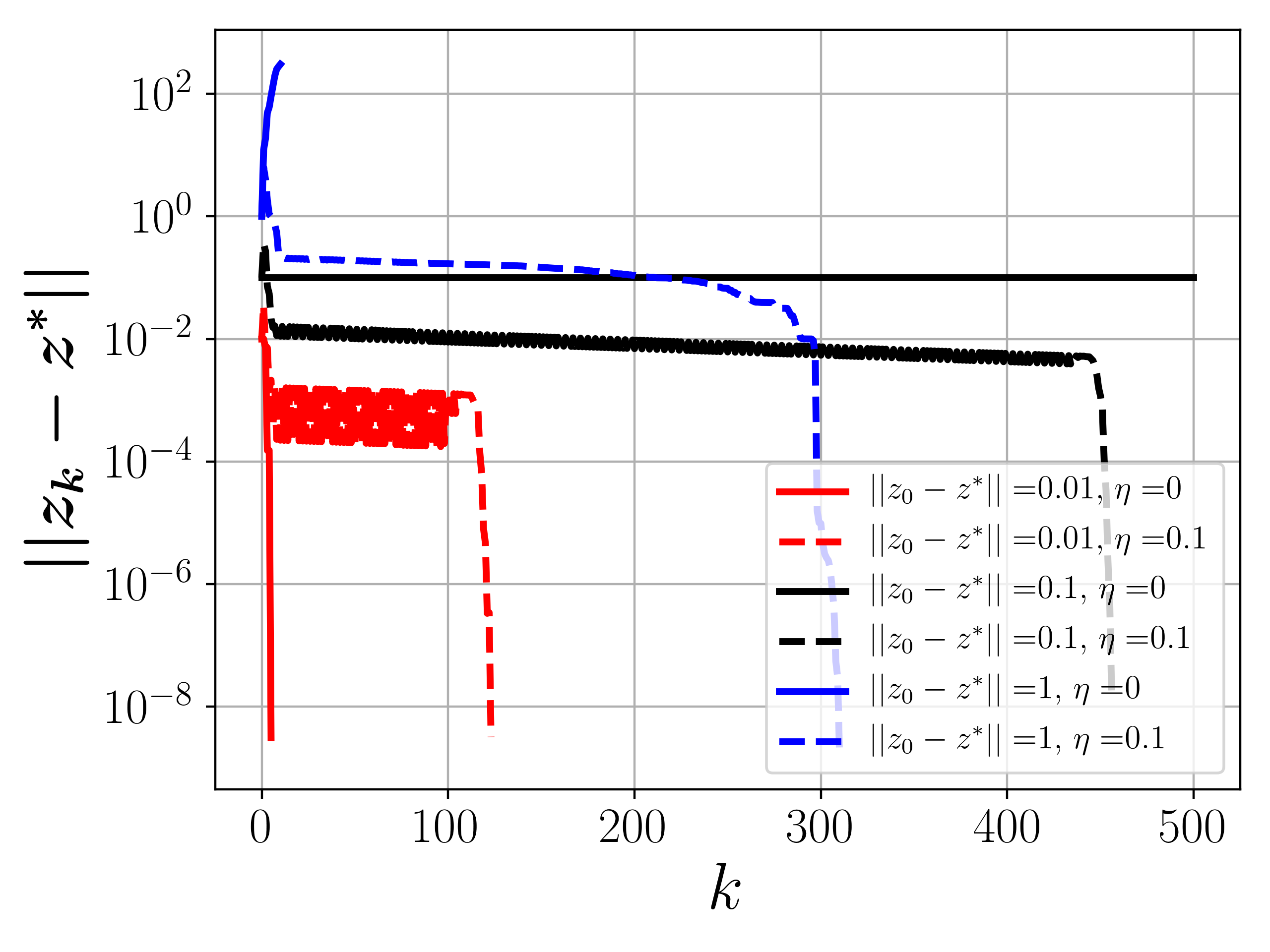}}\label{fig:gan_plot}
\caption{Dirac GAN. Other methods similarly converge to saturation region (leftmost). 
}
\label{fig:dirac}
\end{center}
    \vspace{-0.5cm}
\end{figure}

Generative Adversarial Networks became recently one of the most popular applications of the minimax approach \cite{goodfellow2014generative}.
We test the proposed method in deterministic setting of the Dirac GAN scheme proposed in \cite{mescheder2018training} by expanding the dimensions of the problem such that 
\begin{equation} \label{eq:dirac_gan}
\begin{aligned}
f(x,y) = \phi(-x^{T}y)+\phi(y^{T}c),
\end{aligned}
\end{equation}
for some scalar function $\phi$. Here, $c$ denotes the high dimensional Dirac distribution value that we sample from the Normal distribution. The primal and dual variables represent the data generator and discriminator, respectively.
Figure \ref{fig:dirac} depicts the distance to optimum (leftmost), the gradient norm of the generator and discriminator (second from left), the influence of the subspace dimension (third from left), and the influence of the proximal factor on convergence for different initialization (rightmost).
Therein, we use the common sigmoid cross-entropy loss $\phi(t)=-\ln(1+e^{-t})$ \cite{goodfellow2014generative}.
The dimension of the original optimization problem is set to $M=N=1000$.
Since the objective is concave-concave, the competing methods fail to converge to saddle-point, and diverge to saturation regions \cite{mescheder2018training}. The proximal operator prevents the method from diverging, and allows faster convergence to optimum as depicted in the rightmost plot.
 

\section{Conclusions}
In this paper we introduced a sequential subspace optimization approach to saddle-point problems. 
We improve convergence of first-order methods via efficient secondary subspace optimization.
We evaluated the proposed framework on several saddle-point problems, demonstrating its efficacy and superior performance relative to popular optimization techniques.
Further theoretical investigation of the influence of the subspace directions and dimensions may provide better understanding and enable development of both faster and more efficient saddle-point optimization methods.

\bibliographystyle{plain}
\bibliography{references}

\clearpage
\newpage
\appendix
\section{Proof of Theorem 3.1}
Let us first consider the bi-linear setting $f(x,y)=x^{T}Cy$, where $C$ is a full rank-matrix. We first show that the gradient method is diverging in the above case. The gradient method update is defined as 
\begin{equation}
\begin{aligned}
                \begin{pmatrix}
                  x_{k+1}  \\
                    y_{k+1} 
                \end{pmatrix}
                &=
                \begin{pmatrix}
                  x_{k} \\
                  y_{k}
                \end{pmatrix}
                +\eta_k
                \begin{pmatrix}
                  -\nabla_x f(x_{k},y_{k})\\
                    \nabla_y f(x_{k},y_{k})
                \end{pmatrix}
                \\
                &=
                \begin{pmatrix}
                  x_{k} \\
                  y_{k}
                \end{pmatrix}
                +\eta_k
                \begin{pmatrix}
                  0 & -C\\
                    C^{T} & 0
                \end{pmatrix}
                \begin{pmatrix}
                  x_{k}\\
                    y_{k}
                \end{pmatrix}
                \\
                &=
                \begin{pmatrix}
                  I & -\eta_k C\\
                    \eta_k C^{T} & I
                \end{pmatrix}
                \begin{pmatrix}
                  x_{k}\\
                    y_{k}
                \end{pmatrix}
                =
                A
                \begin{pmatrix}
                  x_{k}\\
                    y_{k}
                \end{pmatrix}
                .
\end{aligned}
\end{equation}
Thus we have $\forall C$ and $\forall \eta_k$
\begin{equation}
\begin{aligned}
            \begin{Vmatrix}
                \begin{pmatrix}
                  x_{k+1}  \\
                    y_{k+1} 
                \end{pmatrix}
            \end{Vmatrix}^2
                &\geq
            \lambda_{\min}(A^{T}A)
            \begin{Vmatrix}
                \begin{pmatrix}
                  x_{k}\\
                    y_{k}
                \end{pmatrix}
            \end{Vmatrix}^2
            =
            \big(1+\lambda_{\min}(\eta_{k}^{2}CC^{T})\big)
            \begin{Vmatrix}
                \begin{pmatrix}
                  x_{k}\\
                    y_{k}
                \end{pmatrix}
            \end{Vmatrix}^2
            .
\end{aligned}
\end{equation}
We now proceed to the proof of convergence of the one dimensional subspace method, assuming $C=I$.
Optimal solution of the subspace optimization satisfies
\begin{equation}
\begin{aligned}
 &\left\{
 \begin{array}{ll}
	    P_k^{T}\nabla_x f(x_{k}+P_k\alpha,y_{k}+Q_k\beta)=0\\
    	Q_k^{T}\nabla_y f(x_{k}+P_k\alpha,y_{k}+Q_k\beta)=0
 \end{array}
 \right.
	 \\ 
	 \Leftrightarrow 
&\left\{
 \begin{array}{ll}
	    \alpha=-(Q_{k}^{T}C^{T}P_{k})^{-1}Q_{k}^{T}C^{T}x_{k}\\
	    \beta=-(P_{k}^{T}CQ_{k})^{-1}P_{k}^{T}Cy_{k}
 \end{array}
 \right.
\end{aligned}
\end{equation}
Thereafter, the update of the variable $x$ is written as
\begin{equation}
\begin{aligned}
x_{k+1} &=x_{k} - \eta_k P_k(Q_{k}^{T}C^{T}P_{k})^{-1}Q_{k}^{T}C^{T}x_{k}=x_{k} - \eta_k P_k(Q_{k}^{T}C^{T}P_{k})^{-1}Q_{k}^{T}Q_{k},
\end{aligned}
\end{equation}
with $\eta_k >0$. Then, we can show that
\begin{equation}
\begin{aligned}
\|x_{k+1}\|^2 &=\|x_{k}\|^2 - 2\frac{\eta_k}{Q_{k}^{T}C^{T}P_{k}} \langle x_{k},P_{k}\|Q_k\|^2\rangle+ \frac{\eta_k^2\|P_k\|^2\|Q_{k}\|^4}{(Q_{k}^{T}C^{T}P_{k})^2}\\
        &=\|x_{k}\|^2 +\frac{\|Q_{k}\|^2}{(Q_{k}^{T}C^{T}P_{k})^2}\Big( -2\eta_{k}Q_{k}^{T}C^{T}P_{k}\langle x_{k},P_{k}\rangle+\eta_{k}^{2}\|P_k\|^2\|Q_{k}\|^2 \Big).
\end{aligned}
\end{equation}
Denoting $\delta_k =-2\eta_{k}Q_{k}^{T}C^{T}P_{k}\langle x,P_{k}\rangle+\eta_{k}^{2}\|P_k\|^2\|Q_{k}\|^2 \\
$ and since $C=I$ we have 
\begin{equation}
\begin{aligned}
\delta_k
=&-2\eta_{k}\|Q_{k}^{T}P_{k}\|^2+\eta_{k}^{2}\|P_k\|^2\|Q_{k}\|^2\\
=& -2\eta_{k}f(x_{k},y_{k})^{2}+\eta_{k}^{2}\|\nabla_x f(x_{k},y_{k})\|^2\|\nabla_y f(x_{k},y_{k})\|^2
\end{aligned}
\end{equation}
Thus, $\forall \eta_{k} \in (0,{2f(x_{k},y_{k})^{2}}/{\|\nabla_x f(x_{k},y_{k})\|^{2}\|\nabla_y f(x_{k},y_{k})\|^{2}})$ we have $\|x_{k+1}\|^2 < \|x_{k}\|^2$.
By following similar arguments, we get $\|y_{k+1}\|^2< \|y_{k}\|^2$.\\
\QEDB


\section{Proof of Theorem 3.2}
We are looking for $\eta>0$ such that $z_{k+1}=z_k +\eta d_{k}$ is a better stationary point than $z_{k}$, i.e. $\|\nabla f(z_{k+1})\| \leq \|\nabla f(z_k)\|$. Here $d_{k}$ designate the anti-gradient $d_x$ and gradient direction $d_y$ according to the primal and dual variable respectively.
From first order expansion we have 
\begin{equation}
\begin{aligned}
f(z_{k+1}) &= f(z_k) + \eta \langle\nabla f(z_k),d_k\rangle + o(\eta^{2}\|d_{k}\|)\\
\iff \nabla{f}(z_{k+1}) &= \nabla f(z_k) + \eta \nabla^{2}f(z_k)d_k
\end{aligned}
\end{equation}
Thus we have
\begin{equation}
\begin{aligned}
\|\nabla f(z_{k+1})\|^2 &= \|\nabla f(z_k)\|^2 + 2\eta \nabla f(z_k)^{T}\nabla^{2} f(z_k)d_k 
+ \eta^{2}d_k^{T}\nabla^2 f(z_k)^{T}\nabla^2 f(z_k) d_k\\
\underset{\eta\rightarrow 0}{\Rightarrow} \|\nabla f(z_{k+1})\|^2 & = \|\nabla f(z_k)\|^2 + 2\eta \nabla f(z_k)^{T}\nabla^{2} f(z_k)d_k,
\end{aligned}
\end{equation}
Since we have 
\begin{equation}
\begin{aligned}
\nabla f(z_k)^{T}\nabla^2 f(z_{k})d_k &= -d_{x}^{T}\nabla_{xx}f(z_k)d_{x} + d_{y}^{T}\nabla_{yy}f(z_k)d_{y} - d_{x}^{T}\nabla_{xy}f(z_k)d_{y} + d_{y}^{T}\nabla_{yx}f(z_k)d_{x} \\
&= -d_{x}^{T}\nabla_{xx}f(z_k)d_{x} + d_{y}^{T}\nabla_{yy}f(z_k)d_{y}<0,
\end{aligned}
\end{equation}
where the last inequality arises from the positive/negative definiteness of the second order partial derivatives in $B_{2}((x^*,y^*),r)$.
\QEDB

Notice the line search procedure cannot diverge for non-strongly convex-concave problems where the block diagonal Hessian can vanish.

\end{document}